\documentclass[12pt]{article}

\usepackage{amsmath,amssymb,amsthm,amscd,a4wide,upref}

\usepackage{hyperref}
\hypersetup{colorlinks,citecolor=blue,filecolor=black,linkcolor=blue,urlcolor=blue}
\usepackage{enumerate}
\usepackage{mathtools}

\begin{document}


\newcommand{\ad}{{\rm ad}}
\newcommand{\cri}{{\rm cri}}
\newcommand{\End}{{\rm{End}\ts}}
\newcommand{\Rep}{{\rm{Rep}\ts}}
\newcommand{\Hom}{{\rm{Hom}}}
\newcommand{\Mat}{{\rm{Mat}}}
\newcommand{\ch}{{\rm{ch}\ts}}
\newcommand{\chara}{{\rm{char}\ts}}
\newcommand{\diag}{{\rm diag}}
\newcommand{\non}{\nonumber}
\newcommand{\wt}{\widetilde}
\newcommand{\wh}{\widehat}
\newcommand{\ot}{\otimes}
\newcommand{\la}{\lambda}
\newcommand{\La}{\Lambda}
\newcommand{\De}{\Delta}
\newcommand{\al}{\alpha}
\newcommand{\be}{\beta}
\newcommand{\ga}{\gamma}
\newcommand{\Ga}{\Gamma}
\newcommand{\ep}{\epsilon}
\newcommand{\ka}{\kappa}
\newcommand{\vk}{\varkappa}
\newcommand{\si}{\sigma}
\newcommand{\vs}{\varsigma}
\newcommand{\vp}{\varphi}
\newcommand{\de}{\delta}
\newcommand{\ze}{\zeta}
\newcommand{\om}{\omega}
\newcommand{\Om}{\Omega}
\newcommand{\ee}{\epsilon^{}}
\newcommand{\su}{s^{}}
\newcommand{\hra}{\hookrightarrow}
\newcommand{\ve}{\varepsilon}
\newcommand{\ts}{\,}
\newcommand{\vac}{\mathbf{1}}
\newcommand{\di}{\partial}
\newcommand{\qin}{q^{-1}}
\newcommand{\tss}{\hspace{1pt}}
\newcommand{\Sr}{ {\rm S}}
\newcommand{\U}{ {\rm U}}
\newcommand{\BL}{ {\overline L}}
\newcommand{\BE}{ {\overline E}}
\newcommand{\BP}{ {\overline P}}
\newcommand{\AAb}{\mathbb{A}\tss}
\newcommand{\CC}{\mathbb{C}\tss}
\newcommand{\KK}{\mathbb{K}\tss}
\newcommand{\QQ}{\mathbb{Q}\tss}
\newcommand{\SSb}{\mathbb{S}\tss}
\newcommand{\TT}{\mathbb{T}\tss}
\newcommand{\ZZ}{\mathbb{Z}\tss}
\newcommand{\DY}{ {\rm DY}}
\newcommand{\X}{ {\rm X}}
\newcommand{\Y}{ {\rm Y}}
\newcommand{\Z}{{\rm Z}}
\newcommand{\Ac}{\mathcal{A}}
\newcommand{\Lc}{\mathcal{L}}
\newcommand{\Mc}{\mathcal{M}}
\newcommand{\Pc}{\mathcal{P}}
\newcommand{\Qc}{\mathcal{Q}}
\newcommand{\Rc}{\mathcal{R}}
\newcommand{\Sc}{\mathcal{S}}
\newcommand{\Tc}{\mathcal{T}}
\newcommand{\Bc}{\mathcal{B}}
\newcommand{\Ec}{\mathcal{E}}
\newcommand{\Fc}{\mathcal{F}}
\newcommand{\Gc}{\mathcal{G}}
\newcommand{\Hc}{\mathcal{H}}
\newcommand{\Uc}{\mathcal{U}}
\newcommand{\Vc}{\mathcal{V}}
\newcommand{\Wc}{\mathcal{W}}
\newcommand{\Yc}{\mathcal{Y}}
\newcommand{\Ar}{{\rm A}}
\newcommand{\Br}{{\rm B}}
\newcommand{\Ir}{{\rm I}}
\newcommand{\Fr}{{\rm F}}
\newcommand{\Jr}{{\rm J}}
\newcommand{\Or}{{\rm O}}
\newcommand{\GL}{{\rm GL}}
\newcommand{\Spr}{{\rm Sp}}
\newcommand{\Rr}{{\rm R}}
\newcommand{\Zr}{{\rm Z}}
\newcommand{\gl}{\mathfrak{gl}}
\newcommand{\middd}{{\rm mid}}
\newcommand{\ev}{{\rm ev}}
\newcommand{\Pf}{{\rm Pf}}
\newcommand{\Norm}{{\rm Norm\tss}}
\newcommand{\oa}{\mathfrak{o}}
\newcommand{\spa}{\mathfrak{sp}}
\newcommand{\osp}{\mathfrak{osp}}
\newcommand{\f}{\mathfrak{f}}
\newcommand{\g}{\mathfrak{g}}
\newcommand{\h}{\mathfrak h}
\newcommand{\n}{\mathfrak n}
\newcommand{\z}{\mathfrak{z}}
\newcommand{\Zgot}{\mathfrak{Z}}
\newcommand{\p}{\mathfrak{p}}
\newcommand{\sll}{\mathfrak{sl}}
\newcommand{\agot}{\mathfrak{a}}
\newcommand{\qdet}{ {\rm qdet}\ts}
\newcommand{\Ber}{ {\rm Ber}\ts}
\newcommand{\HC}{ {\mathcal HC}}
\newcommand{\cdet}{{\rm cdet}}
\newcommand{\rdet}{{\rm rdet}}
\newcommand{\tr}{ {\rm tr}}
\newcommand{\gr}{ {\rm gr}\ts}
\newcommand{\str}{ {\rm str}}
\newcommand{\loc}{{\rm loc}}
\newcommand{\Gr}{{\rm G}}
\newcommand{\sgn}{ {\rm sgn}\ts}
\newcommand{\sign}{{\rm sgn}}
\newcommand{\ba}{\bar{a}}
\newcommand{\bb}{\bar{b}}
\newcommand{\bi}{\bar{\imath}}
\newcommand{\bj}{\bar{\jmath}}
\newcommand{\bk}{\bar{k}}
\newcommand{\bl}{\bar{l}}
\newcommand{\hb}{\mathbf{h}}
\newcommand{\Sym}{\mathfrak S}
\newcommand{\fand}{\quad\text{and}\quad}
\newcommand{\Fand}{\qquad\text{and}\qquad}
\newcommand{\For}{\qquad\text{or}\qquad}
\newcommand{\OR}{\qquad\text{or}\qquad}
\newcommand{\grpr}{{\rm gr}^{\tss\prime}\ts}
\newcommand{\degpr}{{\rm deg}^{\tss\prime}\tss}

\renewcommand{\theequation}{\arabic{section}.\arabic{equation}}

\newtheorem{thm}{Theorem}[section]
\newtheorem{lem}[thm]{Lemma}
\newtheorem{prop}[thm]{Proposition}
\newtheorem{cor}[thm]{Corollary}
\newtheorem{conj}[thm]{Conjecture}
\newtheorem*{mthm}{Main Theorem}
\newtheorem*{mthma}{Theorem A}
\newtheorem*{mthmb}{Theorem B}
\newtheorem*{mthmc}{Theorem C}
\newtheorem*{mthmd}{Theorem D}

\theoremstyle{definition}
\newtheorem{defin}[thm]{Definition}

\theoremstyle{remark}
\newtheorem{remark}[thm]{Remark}
\newtheorem{example}[thm]{Example}

\newcommand{\bth}{\begin{thm}}
\renewcommand{\eth}{\end{thm}}
\newcommand{\bpr}{\begin{prop}}
\newcommand{\epr}{\end{prop}}
\newcommand{\ble}{\begin{lem}}
\newcommand{\ele}{\end{lem}}
\newcommand{\bco}{\begin{cor}}
\newcommand{\eco}{\end{cor}}
\newcommand{\bde}{\begin{defin}}
\newcommand{\ede}{\end{defin}}
\newcommand{\bex}{\begin{example}}
\newcommand{\eex}{\end{example}}
\newcommand{\bre}{\begin{remark}}
\newcommand{\ere}{\end{remark}}
\newcommand{\bcj}{\begin{conj}}
\newcommand{\ecj}{\end{conj}}

\newcommand{\bal}{\begin{aligned}}
\newcommand{\eal}{\end{aligned}}
\newcommand{\beq}{\begin{equation}}
\newcommand{\eeq}{\end{equation}}
\newcommand{\ben}{\begin{equation*}}
\newcommand{\een}{\end{equation*}}

\newcommand{\bpf}{\begin{proof}}
\newcommand{\epf}{\end{proof}}

\def\beql#1{\begin{equation}\label{#1}}


\newcommand{\Res}{\mathop{\mathrm{Res}}}

\title{\Large\bf Center of the quantum affine vertex algebra in type $A$}

\author{{Naihuan Jing,\quad Slaven Ko\v{z}i\'{c},\quad Alexander Molev\quad and\quad Fan Yang}}

\date{} 
\maketitle

\vspace{4 mm}

\begin{abstract}
We consider the quantum vertex algebra associated with the double Yangian
in type $A$ as defined by Etingof and Kazhdan. We show that its
center is a commutative associative algebra and construct
algebraically independent families of topological generators of the center
at the critical level.

\end{abstract}

\vspace{5 mm}


%

\section{Introduction}
\label{sec:int}
\setcounter{equation}{0}

Let $\g$ be a simple Lie algebra over $\CC$
and let $\wh\g$ be the corresponding
affine Kac--Moody algebra.
The vacuum module $V_{\ka}(\g)$ at the level $\ka\in\CC$ over $\wh\g$
has a vertex algebra structure; see, e.g., books by E.~Frenkel and
D.~Ben-Zvi~\cite{fb:va},
I.~Frenkel, J.~Lepowsky and A.~Meurman~\cite{flm:vo} and V.~Kac~\cite{k:va}.
The center of any vertex algebra is a commutative associative algebra.
Unless the level $\ka$ is critical, the center of
the affine vertex algebra $V_{\ka}(\g)$ is trivial (coincides with $\CC$).
By a theorem of B.~Feigin
and E.~Frenkel~\cite{ff:ak}, the center at the
critical level $\z(\wh\g)$
is an algebra of polynomials in infinitely many variables.
Moreover, the algebra $\z(\wh\g)$ is canonically isomorphic to
the algebra of functions on a space of opers; see E.~Frenkel~\cite[Ch.~4]{f:lc} for
a detailed exposition.

Explicit formulas for generators of the Feigin--Frenkel center $\z(\wh\g)$
were given in \cite{cm:ho} and \cite{ct:qs} for type $A$ (see also \cite{mr:mm}),
in \cite{m:ff} for types $B$, $C$ and $D$; and in \cite{mrr:ss} for type $G_2$.
Due to general results of \cite{ffr:gm}, \cite{fft:gm} and \cite{r:si},
these formulas lead to explicit constructions
of commutative subalgebras of the universal enveloping algebra $\U(\g)$
and to explicit higher
order Hamiltonians and their eigenvalues on the Bethe vectors
in the Gaudin model associated with $\g$; see also \cite{fm:qs}, \cite{mm:eb}.

A general definition of quantum vertex algebra was given by P.~Etingof and D.~Kazhdan~\cite{ek:ql5}.
In accordance with \cite{ek:ql5},
a {\em quantum affine vertex algebra} can be associated with a rational,
trigonometric or elliptic $R$-matrix.
In particular, a suitably normalized Yang $R$-matrix gives rise to a
quantum vertex algebra structure on the vacuum module $\Vc_c(\gl_N)$
at the level $c\in\CC$ over the double Yangian
$\DY(\gl_N)$ of type $A$.

In this paper we
introduce the {\em center} $\z(V)$ of an arbitrary quantum vertex algebra $V$
and describe its general properties. We show that the center is an
$\Sc$-{\em commutative} associative
algebra; see \eqref{sloc} for the definition.
Our main focus will be on the center
$\z\big(\Vc_c(\gl_N)\big)$
of the quantum affine vertex algebra $\Vc_c(\gl_N)$. The vacuum
module is isomorphic to the $h$-adically completed dual Yangian
$\Y^+(\gl_N)$, as a vector space, and we prove that the center can be identified with
a {\em commutative} subalgebra of the $h$-adically completed $\Y^+(\gl_N)$.
This subalgebra is invariant under
a derivation $D$, the {\em translation operator}, arising from
the quantum vertex algebra structure on the vacuum module.

We show that the center at the critical level $c=-N$
possesses large families of algebraically independent topological generators so that
a quantum analogue of the Feigin--Frenkel theorem holds.
Moreover, unlike the center of the affine vertex algebra $V_{-N}(\gl_N)$,
it turns out to be possible to produce such families parameterized by arbitrary
partitions with at most $N$ parts. The construction depends of the {\em fusion procedure}
originated in the work of A.~Jucys~\cite{j:yo}
for the symmetric group providing factorized $R$-matrix formulas for all primitive
idempotents. These families thus generalize
to the context of quantum vertex algebras the quantum immanants
of A.~Okounkov~\cite{o:qi} which form a basis of the center of the universal
enveloping algebra $\U(\gl_N)$.

By taking a classical limit we recover explicit generators of the center of the affine
vertex algebra $V_{-N}(\gl_N)$; cf. \cite{cm:ho}, \cite{ct:qs} and \cite{t:qg}.
In principle, this approach is also
applicable to construct generators of the Feigin--Frenkel center $\z(\wh\g)$
for an arbitrary simple Lie algebra $\g$. A required ingredient is a fusion procedure
providing $R$-matrix formulas for idempotents in appropriate centralizer algebras.
This is already in place for the types $B$, $C$ and $D$ so that the construction of \cite{m:ff}
can be reproduced in this way.

We also give a construction of
central elements of a completed
double Yangian at the critical level prompted by the quantum vertex algebra
structure. They are used to show that the center $\z\big(\Vc_c(\gl_N)\big)$
is commutative.
If the level is not critical, then the center
is trivial in the sense that
its generators are elements associated with the center of the Lie algebra $\gl_N$.
They are found as the coefficients of the quantum determinant of the generator matrix
of the dual Yangian.

Our arguments are based on explicit constructions of elements of the center
of the quantum affine vertex algebra
and rely on the $R$-matrix calculations used in \cite{fjmr:hs}
to produce explicit Sugawara operators for the quantum affine algebra
in type $A$ at the critical level.

The research reported in this paper
was supported by the South China University of Technology
in Guangzhou.
This work was finalized during the third author's visit
to the University.
He would like to thank the Center of Quantum Algebra
and the School of Mathematical Sciences
for the warm hospitality during his visit.

\section{Vacuum module for the double Yangian}
\label{sec:dy}
\setcounter{equation}{0}

We recall and reproduce some basic properties of the double Yangian
for $\gl_N$. Our definitions follow
Etingof and Kazhdan~\cite{ek:ql4}, \cite{ek:ql5} and Iohara~\cite{i:br}, where a centrally extended
double Yangian over the ring $\CC[[h]]$ was considered.
To simplify our formulas, we first define this algebra over $\CC$
(formally putting $h=-1$ in the notation of \cite{i:br}), although this will require a certain
completion; cf. Nazarov~\cite{n:yq}. We will return to the closely related
definition of the double Yangian over $\CC[[h]]$ to study the associated structure
of quantum vertex algebra in Sec.~\ref{sec:qava}.

\subsection{Yangian and dual Yangian for $\gl_N$}
\label{subsec:yagln}

The {\it Yangian\/}
$\Y(\gl_N)$
is the associative algebra with generators
$t_{ij}^{(r)}$, where $1\leqslant i,j\leqslant N$ and $r=1,2,\dots$
and the defining relations
\beql{defequiv}
[t^{(r)}_{ij}, t^{(s)}_{kl}] =\sum_{a=1}^{\min\{r,s\}}
\Big(t^{(a-1)}_{kj} t^{(r+s-a)}_{il}-t^{(r+s-a)}_{kj} t^{(a-1)}_{il}\Big),
\eeq
where $t^{(0)}_{ij}=\de_{ij}$. In terms of the formal series
\ben
t_{ij}(u)=\de_{ij}+\sum_{r=1}^{\infty}t_{ij}^{(r)}\ts u^{-r}
\in\Y(\gl_N)[[u^{-1}]]
\een
the defining relations can be written as
\beql{defrel}
(u-v)\ts [t_{ij}(u),t_{kl}(v)]=
t_{kj}(u)\ts t_{il}(v)-t_{kj}(v)\ts t_{il}(u).
\eeq
They admit the following matrix form.
Consider the {\it Yang $R$-matrix\/} $R(u)$, which
is a rational function in a complex parameter $u$
with values in the tensor product algebra
$\End\CC^N\ot\End\CC^N$ defined by
\beql{yangr}
R(u)=1-P\ts u^{-1},
\eeq
where $P$ is the permutation operator in $\CC^N\ot\CC^N$.
Then \eqref{defequiv} is equivalent to the {\em RTT relation}
\beql{RTT}
R(u-v)\ts T_1(u)\ts T_2(v)=T_2(v)\ts T_1(u)\ts R(u-v),
\eeq
where
\beql{tu}
T(u)=\sum_{i,j=1}^N e_{ij}\ot t_{ij}(u)
\in \End\CC^N\ot \Y(\gl_N)[[u^{-1}]]
\eeq
and the $e_{ij}$ are the matrix units.
We use a subscript to indicate a copy of the matrix of the form \eqref{tu} in the multiple
tensor product algebra
\beql{multtpr}
\underbrace{\End\CC^N\ot\dots\ot\End\CC^N}_m\ot\Y(\gl_N)[[u^{-1}]]
\eeq
so that
\beql{la}
T_a(u)=\sum_{i,j=1}^N 1^{\ot (a-1)}\ot e_{ij}\ot 1^{\ot (m-a)}\ot t_{ij}(u).
\eeq
We take $m=2$ for the defining relations \eqref{RTT}.

This notation for elements of algebras of the form \eqref{multtpr} will be extended
as follows. For an element
\ben
C=\sum_{i,j,r,s=1}^N c^{}_{ijrs}\ts e_{ij}\ot e_{rs}\in
\End \CC^N\ot\End \CC^N,
\een
and any two indices $a,b\in\{1,\dots,m\}$ such that $a\ne b$,
we denote by $C_{a\tss b}$ the element of the algebra $(\End\CC^N)^{\ot m}$ with $m\geqslant 2$
given by
\beql{cab}
C_{a\tss b}=\sum_{i,j,r,s=1}^N c^{}_{ijrs}\ts (e_{ij})_a\tss (e_{rs})_b,
\qquad
(e_{ij})_a=1^{\ot(a-1)}\ot e_{ij}\ot 1^{\ot(m-a)}.
\eeq
We regard the matrix transposition as the linear map
\ben
t:\End\CC^N\to\End\CC^N,\qquad e_{ij}\mapsto e_{ji}.
\een
For any $a\in\{1,\dots,m\}$ we will denote by $t_a$ the corresponding
partial transposition on the algebra \eqref{multtpr} which acts as $t$ on the
$a$-th copy of $\End \CC^N$ and as the identity map on all the other tensor factors.

The algebra $\Y(\gl_N)$ possesses a natural ascending filtration
defined by $\deg t_{ij}^{(r)}=r-1$
for all $r\geqslant1$. Denote by
$\gr\Y(\gl_N)$ the associated graded algebra.
We have the isomorphism
$\gr\Y(\gl_N)\cong \U\big(\gl_N[t]\big)$.
The image $\bar t_{ij}^{\ts(r)}$ of the generator $t_{ij}^{(r)}$ in the $(r-1)$-th
component of the graded algebra $\gr\Y(\gl_N)$ corresponds to the element $E_{ij}[r-1]$
of $\U\big(\gl_N[t]\big)$, where the $E_{ij}$ are the standard basis elements of $\gl_N$ and
we use the notation $X[r]=X\tss t^r$ for
$X\in\gl_N$ and any $r\in\ZZ$.

Let $E=[E_{ij}]$ denote the matrix whose $(i,j)$ entry
is the element $E_{ij}$ of $\U(\gl_N)$. For any $a\in\CC$ the mapping
\beql{evalyu}
\ev_a:T(u)\mapsto 1+E\ts (u-a)^{-1},
\eeq
defines a homomorphism $\Y(\gl_N)\to\U(\gl_N)$ known as
the {\it evaluation homomorphism\/}. In terms of generators,
 $\ev_a:t_{ij}^{\ts(r)}\mapsto E_{ij}\ts a^{r-1}$.

For more details on
the origins, structure
and representations of the Yangian see \cite{m:yc}.

The
{\em dual Yangian}
$\Y^+(\gl_N)$ can be defined
as the associative algebra with generators
$t_{ij}^{(-r)}$, where $1\leqslant i,j\leqslant N$ and $r=1,2,\dots$
subject to the defining relations
\beql{defdual}
[t^{(-r)}_{ij}, t^{(-s)}_{kl}] =\de_{kj}\ts t^{(-r-s)}_{il}- \de_{il}\ts t^{(-r-s)}_{kj}
+\sum_{a=1}^{\min\{r,s\}}
\Big(t^{(-r-s+a-1)}_{kj} t^{(-a)}_{il}-t^{(-a)}_{kj} t^{(-r-s+a-1)}_{il}\Big).
\eeq
Combining the generators into the formal power series
\ben
t^+_{ij}(u)=\de_{ij}-\sum_{r=1}^{\infty}t_{ij}^{(-r)}\ts u^{r-1}
\in\Y^+(\gl_N)[[u]]
\een
we can write the defining relations as
\beql{defreldual}
(u-v)\ts [t^+_{ij}(u),t^+_{kl}(v)]=
t^+_{kj}(u)\ts t^+_{il}(v)-t^+_{kj}(v)\ts t^+_{il}(u)
\eeq
which thus take the same form as \eqref{defrel}. So they are equivalent to
\beql{RTTdual}
R(u-v)\ts T^+_1(u)\ts T^+_2(v)=T^+_2(v)\ts T^+_1(u)\ts R(u-v)
\eeq
as in \eqref{RTT}, where we use the Yang $R$-matrix \eqref{yangr} and
\beql{tuplus}
T^+(u)=\sum_{i,j=1}^N e_{ij}\ot t^+_{ij}(u)
\in \End\CC^N\ot \Y^+(\gl_N)[[u]].
\eeq

Consider the ascending filtration on the dual Yangian $\Y^+(\gl_N)$
defined by $\deg t_{ij}^{(-r)}=-r$ for all $r\geqslant1$.
We have the isomorphism for the associated graded algebra
\beql{grisom}
\gr\Y^+(\gl_N)\cong \U\big(t^{-1}\gl_N[t^{-1}]\big).
\eeq
The image $\bar t_{ij}^{\ts(-r)}$ of the generator
$t_{ij}^{(-r)}$ in the $(-r)$-th
component of the graded algebra $\gr\Y^+(\gl_N)$ corresponds to the element $E_{ij}[-r]$
of $\U\big(t^{-1}\gl_N[t^{-1}]\big)$.
The isomorphism relies on the Poincar\'e--Birkhoff--Witt theorem
for $\Y^+(\gl_N)$ which can be proved in a way similar to the Yangian; cf.~\cite[Ch.~1]{m:yc}
and references therein. We will give a more general proof below in the context
of the double Yangian which would imply \eqref{grisom}; see Corollary~\ref{cor:isompbw}.
For any nonzero $a\in\CC$ the mapping
\beql{evalyudu}
\ev_a:T^+(u)\mapsto 1+E\ts (u-a)^{-1},
\eeq
defines the {\it evaluation homomorphism\/}
$\Y^+(\gl_N)\to\U(\gl_N)$. We assume an expansion into a power series in $u$ so that
in terms of generators it takes the form
 $\ev_a:t_{ij}^{\ts(-r)}\mapsto E_{ij}\ts a^{-r}$.

\subsection{Double Yangian for $\gl_N$}
\label{subsec:dougln}

The {\em double Yangian} $\DY(\gl_N)$ for $\gl_N$ is defined as the associative algebra
generated by the central element $C$ and elements
$t_{ij}^{(r)}$ and $t_{ij}^{(-r)}$, where $1\leqslant i,j\leqslant N$ and $r=1,2,\dots$,
subject to the defining relations written in terms of the generator matrices
\eqref{tu} and \eqref{tuplus} as follows; see \cite{ek:ql4}, \cite{ek:ql5} and \cite{i:br}.
They are given by \eqref{RTT}, \eqref{RTTdual}
together with the relation
\beql{RTTdou}
\overline R\big(u-v+C/2\big)\ts T_1(u)\ts T^+_2(v)=T^+_2(v)\ts T_1(u)\ts
\overline R\big(u-v-C/2\big),
\eeq
where
\beql{orm}
\overline R(u)=g(u) \ts R(u)=g(u)\ts\big(1-P\ts u^{-1}\big)
\eeq
and
\beql{gu}
g(u)=1+\sum_{i=1}^{\infty}\ts g_i\ts u^{-i},\qquad g_i\in\CC,
\eeq
is a formal power series in $u^{-1}$ whose coefficients are
uniquely determined by the relation
\beql{recg}
g(u+N)=g(u)\ts (1-u^{-2}).
\eeq
Its first few terms are
\ben
g(u)=1+\frac{1}{N}\ts u^{-1}+\frac{N^2+1}{2\tss N^2}\ts u^{-2}
+\frac{N^4+4\tss N^2+1}{6\tss N^3}\ts u^{-3}+\dots.
\een
The relation \eqref{recg} ensures that the $R$-matrix $\overline R(u)=\overline R_{12}(u)$ possesses
the {\em crossing symmetry} properties
\beql{cs}
\big(\overline R_{12}(u)^{-1}\big)^{t_1} \ts \overline R_{12}\big(u+N\big)^{t_1}=1
\Fand
\big(\overline R_{12}(u)^{-1}\big)^{t_2} \ts \overline R_{12}\big(u+N\big)^{t_2}=1.
\eeq
Moreover, the following
{\em unitarity property} holds
\beql{up}
\overline R_{12}(u)\tss \overline R_{12}(-u)=1.
\eeq
Indeed, replacing $u$ with $-u-N$ in \eqref{recg} we get
\ben
g(-u)=g(-u-N)\ts \big(1-(u+N)^{-2}\big)
\een
and so
\ben
g(u)\tss g(-u)\ts (1-u^{-2})=g(u+N)\tss g(-u-N)\tss \big(1-(u+N)^{-2}\big).
\een
This means that the series on the left hand side is invariant under the
shift $u\mapsto u+N$ which is only possible when
\ben
g(u)\tss g(-u)\ts (1-u^{-2})=1
\een
thus implying \eqref{up}.
The series $g(u)$ can be defined equivalently as a unique
formal power series of the form \eqref{gu} satisfying the relation
\beql{guprod}
g(u)\ts g(u+1)\dots g(u+N-1)=\big(1-u^{-1}\big)^{-1}.
\eeq
To see the equivalence of the definitions, observe that by \eqref{recg},
the series $G(u)$ defined by the left hand side of \eqref{guprod}
satisfies $G(u+1)=G(u)\ts (1-u^{-2})$. However, $G(u)$ is uniquely determined
by this relation and so coincides with the right hand side of \eqref{guprod}.

Given any $c\in\CC$ we will introduce the {\em double Yangian at the level} $c$
as the quotient $\DY_c(\gl_N)$ of $\DY(\gl_N)$
by the ideal generated by $C-c$. In particular, we have the natural
epimorphism
\beql{epicze}
\vp:\DY(\gl_N)\to \DY_0(\gl_N),\qquad C\mapsto 0,\quad t_{ij}^{(r)}\mapsto t_{ij}^{(r)}.
\eeq

Equip the Lie algebra $\gl_N$ with
the invariant symmetric bilinear form given by
\ben
\langle X,\ts Y\rangle=\tr\tss(X\tss Y)-\frac{1}{N}\ts\tr\tss X
\ts\tr\tss Y,\qquad X, Y\in\gl_N.
\een
Consider the corresponding
affine
Kac--Moody algebra $\wh\gl_N=\gl_N[t,t^{-1}]\oplus\tss\CC K$
defined by
the commutation relations
\beql{commrel}
\big[E_{ij}[r],E_{kl}[s\tss]\tss\big]
=\de_{kj}\ts E_{i\tss l}[r+s\tss]
-\de_{i\tss l}\ts E_{kj}[r+s\tss]
+r\tss\de_{r,-s}\ts K\Big(\de_{kj}\tss\de_{i\tss l}
-\frac{\de_{ij}\tss\de_{kl}}{N}\Big),
\eeq
and the element $K$ is central.

Introduce the ascending filtration on the double Yangian $\DY(\gl_N)$ by
\beql{ascfi}
\deg t_{ij}^{(r)}=r-1\Fand \deg t_{ij}^{(-r)}=-r
\eeq
for all $r\geqslant1$; the degree of the central element $C$ is defined to be equal to zero.
Denote by
$\gr\DY(\gl_N)$ the corresponding graded algebra. We will use the notation
$\bar t_{ij}^{\ts(r)}$ and $\bar t_{ij}^{\ts(-r)}$ for the images of the generators
in the respective components of the graded algebra and let $\overline C$ be the image of $C$
in the zeroth component.

\bpr\label{prop:grdy}
The assignments
\beql{assigttc}
E_{ij}[r-1]\mapsto\bar t_{ij}^{\ts(r)},\qquad E_{ij}[-r]\mapsto\bar t_{ij}^{\ts(-r)}
\Fand K\mapsto\overline C
\eeq
with $r\geqslant 1$ define a homomorphism
\beql{homim}
\U(\wh\gl_N)\to\gr\DY(\gl_N).
\eeq
\epr

\bpf
As we pointed out in the previous sections, there are homomorphisms
\ben
\U\big(\gl_N[t]\big)\to\gr\Y(\gl_N)\Fand \U\big(t^{-1}\gl_N[t^{-1}]\big)\to\gr\Y^+(\gl_N)
\een
which are defined by the assignments \eqref{assigttc}.
We will now use the defining relations \eqref{RTTdou} to verify that
the generators $\bar t_{ij}^{\ts(r)}$ and $\bar t_{kl}^{\ts(-s)}$
with $r,s\geqslant 1$ of the graded algebra satisfy
the desired relations in $\U(\wh\gl_N)$.
Introduce re-scaled generators of $\DY(\gl_N)$ by setting
\ben
\wt t_{ij}^{\ts(r)}=h^{r-1}\ts t_{ij}^{(r)}\Fand \wt t_{ij}^{\ts(-r)}=h^{-r}\ts t_{ij}^{\ts(-r)}
\een
for $r\geqslant 1$, where $h$ is a complex-valued parameter.
The relations satisfied by $\bar t_{ij}^{\ts(r)}$ and $\bar t_{kl}^{\ts(-s)}$
in the graded algebra $\gr\DY(\gl_N)$ will be recovered by calculating
the relations between $\wt t_{ij}^{\ts(r)}$ and $\wt t_{kl}^{\ts(-s)}$
and then taking the limit as $h\to 0$.
Set
\ben
\wt t_{ij}(u)=\sum_{r=1}^{\infty}\wt t_{ij}^{\ts(r)}\ts u^{-r}
=\frac{1}{h}\ts\Big(t_{ij}\Big(\frac{u}{h}\Big)-\de_{ij}\Big)
\een
and
\ben
\wt t^+_{kl}(v)=\sum_{s=1}^{\infty}\wt t_{kl}^{\ts(-s)}\ts v^{s-1}
=\frac{1}{h}\ts\Big(\de_{kl}-t^+_{kl}\Big(\frac{v}{h}\Big)\Big).
\een
Write \eqref{RTTdou} in terms of the generating series:
\begin{multline}
\label{ttplus}
g\big(u-v+C/2\big)\ts \Big(t_{ij}(u)\ts t^+_{kl}(v)-\frac{1}{u-v+C/2}\ts
t_{kj}(u)\ts t^+_{il}(v)\Big)\\[0.5em]
{}=g\big(u-v-C/2\big)\ts \Big(t^+_{kl}(v)\ts t_{ij}(u)-\frac{1}{u-v-C/2}\ts
t^+_{kj}(v)\ts t_{il}(u)\Big).
\end{multline}
Note the expansion into a power series in $(u-v)^{-1}$:
\beql{ratg}
\frac{g\big(u-v-C/2\big)}{g\big(u-v+C/2\big)}=1+\frac{C}{N\tss(u-v)^2}+\dots.
\eeq
Now replace $u$ by $u/h$ and $v$ by $v/h$ in \eqref{ttplus} to get the corresponding
relations between the series
$\wt t_{ij}(u)$ and $\wt t^+_{kl}(v)$. We have
\ben
\bal
\big(\de_{ij}+h\tss \wt t_{ij}(u)\big)\big(\de_{kl}-h\tss \wt t^+_{kl}(v)\big)
&-\frac{h}{u-v+h\tss C/2}\ts
\big(\de_{kj}+h\tss \wt t_{kj}(u)\big)\big(\de_{il}-h\tss \wt t^+_{il}(v)\big)\\[0.5em]
{}-\Big(\big(\de_{kl}-h\tss \wt t^+_{kl}(v)\big)\big(\de_{ij}+h\tss \wt t_{ij}(u)\big)
&-\frac{h}{u-v-h\tss C/2}\ts
\big(\de_{kj}-h\tss \wt t^+_{kj}(v)\big)\big(\de_{il}+h\tss \wt t_{il}(u)\big)\Big)\\[0.5em]
{}&\times\Big(1+\frac{h^2\tss C}{N\tss(u-v)^2}+\dots\Big)=0.
\eal
\een
As a power series in $h$, the left hand side is divisible by
$h^2$. Hence, dividing by $h^2$ we get the relation modulo $h$,
\ben
\bal
\big[\wt t_{ij}(u),\wt t^+_{kl}(v)\big]&\equiv\frac{1}{u-v}
\Big(\de_{kj}\tss\big(\ts\wt t_{il}(u)+\wt t^+_{il}(v)\big)
-\de_{il}\tss\big(\ts\wt t_{kj}(u)+\wt t^+_{kj}(v)\big)\Big)\\[0.5em]
{}&+\frac{C}{N\tss(u-v)^2}\Big(N\tss\de_{kj}\tss \de_{il}-\de_{ij}\tss \de_{kl}\Big).
\eal
\een
Thus, taking the coefficients of $u^{-r}v^{s-1}$ with $r,s\geqslant 1$ on both sides, in
the limit $h\to 0$ in the graded algebra we get
\ben
\big[\bar t_{ij}^{\ts(r)},\bar t_{kl}^{\ts(-s)}\tss\big]
=\begin{cases}\de_{kj}\ts \bar t_{i\tss l}^{\ts(r-s)}
-\de_{i\tss l}\ts \bar t_{k\tss j}^{\ts(r-s)}
+(r-1)\ts\de_{r,s+1}\ts \overline C\Big(\de_{kj}\tss\de_{i\tss l}
-\dfrac{\de_{ij}\tss\de_{kl}}{N}\Big)\quad&\text{if}\quad r>s,\\[0.9em]
\de_{kj}\ts \bar t_{i\tss l}^{\ts(r-s-1)}
-\de_{i\tss l}\ts \bar t_{k\tss j}^{\ts(r-s-1)}
+(r-1)\ts\de_{r,s+1}\ts \overline C\Big(\de_{kj}\tss\de_{i\tss l}
-\dfrac{\de_{ij}\tss\de_{kl}}{N}\Big)\quad&\text{if}\quad r\leqslant s.
\end{cases}
\een
Comparing with \eqref{commrel}, we may conclude that the
assignments \eqref{assigttc} define a homomorphism
\eqref{homim}.
\epf

To equip the double Yangian with a Hopf algebra structure, we will need to
use shifts $u\mapsto u+a$ of the variable $u$. They are well-defined for the generator
series $t_{ij}(u)$ but not for $t^+_{ij}(u)$. So we will consider
the completion $\wh\Y^+(\gl_N)$ of the dual Yangian with respect to the
{\em descending} filtration defined
by setting
the degree of $t_{ij}^{(-r)}$ with $r\geqslant 1$
to be equal to $r$. By the defining relations,
the double Yangian $\DY(\gl_N)$ is spanned over $\CC[C]$ by the products
$x\tss y$ with $x\in \Y^+(\gl_N)$ and $y\in \Y(\gl_N)$. This follows by an easy induction
based on the relation obtained by
swapping the indices $i$ and $k$
in \eqref{ttplus} and
solving the system of equations for $t_{ij}(u)\ts t^+_{kl}(v)$ and
$t_{kj}(u)\ts t^+_{il}(v)$.
The {\em extended double Yangian} $\DY^{\circ}(\gl_N)$
can now be
defined as the space of finite $\CC[C]$-linear combinations of all products of the form $x\tss y$
with $x\in\wh\Y^+(\gl_N)$ and $y\in\Y(\gl_N)$ with the multiplication extended by continuity
from the double Yangian.

The Hopf algebra structure on $\DY^{\circ}(\gl_N)$ is defined by
the coproduct
\begin{align*}
&\Delta:t_{ij}(u)\mapsto\sum_{k=1}^{N} t_{ik}\big(u+ C_2 /4\big)
\otimes t_{kj}\big(u- C_1 /4\big),\\
&\Delta:t_{ij}^{+ }(u)\mapsto\sum_{k=1}^{N} t_{ik}^{+ }\big(u- C_2 /4\big)
\otimes t_{kj}^{+ }\big(u+ C_1 /4\big),\\[0.3em]
&\Delta:C\mapsto C\otimes 1+1\otimes C,
\end{align*}
where $C_1=C\otimes 1$ and $C_2 =1\otimes C$; the antipode
\ben
S:T(u)\mapsto T(u)^{-1},\qquad S:T^{+}(u)\mapsto T^{+}(u)^{-1},\qquad S:C\mapsto -C;
\een
and the counit
\ben
\varepsilon:T(u)\mapsto 1,\qquad \varepsilon:T^{+}(u)\mapsto 1,\qquad\varepsilon:C\mapsto 0.
\een

We are now in a position to prove the Poincar\'e--Birkhoff--Witt theorem for the double Yangian.
We will do this for a particular ordering $\prec$ on the set of generators.
Observe that due to the defining relations,
\ben
\big[t^{(r)}_{ij},t^{(s)}_{ij}\big]=0\Fand \big[t^{(-r)}_{ij},t^{(-s)}_{ij}\big]=0
\een
for all $r,s\geqslant 1$. Hence, a total ordering $\prec$ on the series $t_{ij}(u)$
and $t^+_{ij}(u)$ will induce a well-defined total ordering on the generators
(with the central element $C$ included in the ordering in an arbitrary way).
We set $t^+_{ij}(u)\prec t_{kl}(u)$ for all $i,j,k,l$. Furthermore, set
$t^+_{ij}(u)\prec t^+_{kl}(u)$ and $t_{ij}(u)\prec t_{kl}(u)$
if and only if $(i,j)\prec(k,l)$ in the lexicographical order.

\bth\label{thm:pbw}
Any element of the algebra $\DY(\gl_N)$ can be
written uniquely as a linear combination of ordered monomials in the generators.
\eth

\bpf
It follows by an easy induction from the defining relations \eqref{defrel}, \eqref{defreldual}
and \eqref{ttplus} that the ordered monomials span the algebra $\DY(\gl_N)$.

The next step is to demonstrate that the ordered monomials are linearly independent.
We consider the level zero algebra $\DY_0(\gl_N)$ first and follow
an idea used by Etingof and Kazhdan~\cite[Proposition~3.15]{ek:ql3} and
by Nazarov~\cite[Proposition~2.2]{n:yq}.
It is based on the existence of the evaluation modules for $\DY_0(\gl_N)$:
for each nonzero $a\in\CC$
we have the representation defined by
\ben
\pi_a:\DY_0(\gl_N)\to\End\CC^N,\qquad
t^{(r)}_{ij}\mapsto a^{\tss r-1}\tss e_{ij},\quad
t^{(-r)}_{ij}\mapsto a^{\tss -r}\tss e_{ij}
\een
for all $r\geqslant 1$. If there is a nontrivial linear combination
of ordered monomials equal to zero, we employ
the coproduct on $\DY_0(\gl_N)$ to conclude that the image of this linear combination
is zero under any representation $\pi_{a_1}\ot\dots\ot\pi_{a_l}$ with nonzero
parameters $a_i$. This leads to a contradiction
exactly as in \cite{n:yq} by considering the top
degree components of all monomials with respect to the filtration defined
by \eqref{ascfi} and by employing associated
evaluation modules over $\U\big(\gl_N\tss[t,t^{-1}]\big)$
as implied by Proposition~\ref{prop:grdy}.

To show that ordered monomials are linearly independent in
$\DY(\gl_N)$, observe that $C\ne 0$ due to the existence of
the level $1$ representations. Here we rely on the work by Iohara~\cite{i:br}
providing such representations in terms of the Drinfeld presentation of $\DY(\gl_N)$
to be written in terms of the $RTT$ presentation via the Ding--Frenkel
isomorphism. Now prove by the induction on $k\geqslant 1$ that the powers
$1,C,\dots,C^k$ are linearly independent. Suppose that
\ben
d_k\tss C^k+\dots+d_1\tss C+d_0=0,\qquad d_i\in\CC,\quad d_k\ne 0.
\een
By applying the homomorphism $\vp$ defined in
\eqref{epicze} we find that $d_0=0$. If $k=1$ then
this makes a contradiction since $C\ne 0$. Now suppose that $k\geqslant 2$
and apply the coproduct map $\De$ to get
\ben
d_k\tss (C\ot 1+1\ot C)^k+\dots+d_1\tss (C\ot 1+1\ot C)=0.
\een
This simplifies to
\ben
d_k\tss (k\tss C^{k-1}\ot C+\dots+k\tss C\ot C^{k-1})+\dots+2\tss d_2\tss C\ot C=0.
\een
However,
this is impossible since the powers $1,C,\dots,C^{k-1}$ are linearly independent
by the induction hypothesis.
Now suppose that a linear combination of ordered monomials is zero,
\beql{lcoac}
A(C)+
\sum A_{i_1j_1\dots\ts i_pj_p}^{r_1\dots\ts r_p}(C)\tss t^{(r_1)}_{i_1j_1}\dots t^{(r_p)}_{i_pj_p}=0,
\eeq
where the summation is over a finite nonempty set of indices,
$A(C)$ is a polynomial in $C$, and
the coefficients $A_{i_1j_1\dots\ts i_pj_p}^{r_1\dots\ts r_p}(C)$ are nonzero polynomials in $C$.
Regarding $\DY(\gl_N)$ as a subalgebra of $\DY^{\circ}(\gl_N)$,
apply the homomorphism $\psi=(\text{\rm id\ts}\ot\vp)\ts\De$ to its elements.
The action on polynomials in $C$
is given by
\ben
\psi:B(C)\mapsto B(C)\ot 1,
\een
whereas the images of the generators $t^{(r)}_{ij}$ under $\psi$ are found from the expansions
\ben
t_{ij}(u)\mapsto\sum_{k=1}^{N} t_{ik}(u)
\ot t_{kj}\big(u- C_1 /4\big)\Fand
t_{ij}^{+ }(u)\mapsto\sum_{k=1}^{N} t_{ik}^{+ }(u)
\ot t_{kj}^{+ }\big(u+ C_1 /4\big).
\een
Let $p_0$ be the maximum length of the monomials occurring in
the linear combination in \eqref{lcoac} and let $r_0$ be the maximum degree among
the monomials of length $p_0$. Now apply the homomorphism $\psi$
to the left hand side of \eqref{lcoac} and
use the defining relations of the double Yangian to write the image as a
(possibly infinite) linear
combination of products of the form $x\ot y$, where $x$
and $y$ are ordered monomials in the generators. The defining relations
and coproduct formulas imply that the part of this linear combination
containing the monomials $y$ of length $p_0$ and degree $r_0$ has the form
\beql{linco}
\sum A_{i_1j_1\dots\ts i_pj_p}^{r_1\dots\ts r_p}(C)\ot
 t^{(r_1)}_{i_1j_1}\dots t^{(r_p)}_{i_pj_p},
\eeq
where $p=p_0$ and the sum of the degrees of the generators is equal to $r_0$.
On the other hand, the ordered monomials $t^{(r_1)}_{i_1j_1}\dots t^{(r_p)}_{i_pj_p}$
are linearly independent in $\DY(\gl_N)$ over $\CC$, as follows from
the Poincar\'e--Birkhoff--Witt theorem for the level zero algebra $\DY_0(\gl_N)$
by the application of the homomorphism \eqref{epicze}.
This implies that the coefficients $A_{i_1j_1\dots\ts i_pj_p}^{r_1\dots\ts r_p}(C)$ in
\eqref{linco} must be zero,
thus making a contradiction with the assumptions in \eqref{lcoac}.
Therefore, all ordered monomials
are linearly independent.
\epf

\bco\label{cor:isompbw}
The homomorphism \eqref{homim} is injective and so it defines
an isomorphism
\beql{homimi}
\U(\wh\gl_N)\cong\gr\DY(\gl_N).
\eeq
\eco

\bpf
This is immediate from Theorem~\ref{thm:pbw} and the Poincar\'e--Birkhoff--Witt theorem
for the algebra $\U(\wh\gl_N)$.
\epf

\subsection{Invariants of the extended vacuum module}
\label{subsec:ivmdy}

Theorem~\ref{thm:pbw} implies the vector space decomposition for the extended
double
Yangian as a $\CC[C]$-module,
\beql{dydec}
\DY^{\circ}(\gl_N)\cong \wh\Y^+(\gl_N)\ot \Y(\gl_N).
\eeq
Introduce
the {\em extended vacuum module $\wh\Vc_c(\gl_N)$ at the level $c$}
as the quotient of the algebra $\DY^{\circ}(\gl_N)$ by the left ideal
generated by $C-c$ and all elements $t_{ij}^{\ts(r)}$ with $r\geqslant 1$.
We let $\vac$ denote the image of $1$ in the quotient.
As a vector space, $\wh\Vc_c(\gl_N)$ is isomorphic to the completed dual Yangian $\wh\Y^+(\gl_N)$
due to the decomposition \eqref{dydec}.

Now assume that the level is critical, $c=-N$, and set $\wh\Vc_{\text{cri}}=\wh\Vc_{-N}(\gl_N)$.
Introduce the subspace of $\Y(\gl_N)$-invariants by
\beql{ffvacy}
\z(\wh\Vc_{\text{cri}})=
\{v\in \wh\Vc_{\text{cri}}\ |\ t_{ij}(u)\tss v=\de_{ij}\tss v\},
\eeq
so that any element of $\z(\wh\Vc_{\text{cri}})$ is annihilated by all operators
$t_{ij}^{\ts(r)}$ with $r\geqslant 1$. We will discuss the structure of
the space $\z(\wh\Vc_{\text{cri}})$ below in Sec.~\ref{subsec:cava} in the context of
quantum vertex algebra structure
on $\Vc_{-N}(\gl_N)$. In particular, we will see that
$\z(\wh\Vc_{\text{cri}})$ is a commutative associative algebra which can be
identified with a subalgebra of the completed dual Yangian $\wh\Y^+(\gl_N)$.
Our goal in this section is to construct some families of elements
of $\z(\wh\Vc_{\text{cri}})$.

We will work with the tensor product algebra
\beql{tenprkeyadual}
\underbrace{\End\CC^{N}\ot\dots\ot\End\CC^{N}}_m{}\ot\wh\Y^+(\gl_N)
\eeq
and introduce the rational function in variables $u_1,\dots,u_m$
with values in \eqref{tenprkeyadual} (with the identity component in $\wh\Y^+(\gl_N)$)
by
\beql{rtrain}
R(u_1,\dots,u_m)=\prod_{1\leqslant a<b\leqslant m} R_{a\tss b}(u_a-u_b),
\eeq
where the product is taken in
the lexicographical order on the set of pairs $(a,b)$. We point out
the identity
\beql{rtraintt}
R(u_1,\dots,u_m)\tss T^+_1(u_1)\dots T^+_m(u_m)
=T^+_m(u_m)\dots T^+_1(u_1) R(u_1,\dots,u_m),
\eeq
implied by a repeated application of \eqref{RTTdual}.

Suppose that $\mu$ is a Young diagram with $m$ boxes whose length does not exceed $N$.
For a standard $\mu$-tableau $\Uc$ with entries in $\{1,\dots,m\}$
introduce the contents
$c_a=c_a(\Uc)$ for $a=1,\dots,m$ so that
$c_a=j-i$ if $a$ occupies the box $(i,j)$ in $\Uc$.
Let $e^{}_{\Uc}\in\CC[\Sym_m]$ be the primitive idempotent associated with $\Uc$
through the use of the orthonormal Young bases in the irreducible representations of $\Sym_m$.
The symmetric group $\Sym_m$ acts by permuting the tensor
factors in $(\CC^N)^{\ot\ts m}$. Denote by $\Ec^{}_{\Uc}$ the image
of $e^{}_{\Uc}$ under this action. We will need an expression for $\Ec^{}_{\Uc}$
provided by the fusion procedure originated in \cite{j:yo};
see also \cite[Sec.~6.4]{m:yc} for
more details
and references. By a version of the procedure,
the consecutive evaluations
of the function $R(u_1,\dots,u_m)$ are well-defined
and the result is proportional to $\Ec^{}_{\Uc}$,
\beql{fusion}
R(u_1,\dots,u_m)\big|_{u_1=c_1}
\big|_{u_2=c_2}\dots \big|_{u_m=c_m}=h(\mu)\ts \Ec^{}_{\Uc},
\eeq
where $h(\mu)$ is the product of all hook lengths of the boxes of $\mu$.

Using the tensor product algebra \eqref{tenprkeyadual}, set
\beql{smuuyadual}
\TT^+_{\mu}(u)=\tr^{}_{1,\dots,m}\ts \Ec^{}_{\Uc}\ts T^+_1(u+c_1)\dots T^+_m(u+c_m).
\eeq
This is a power series in $u$ whose coefficients are elements of the completed
dual Yangian $\wh\Y^+(\gl_N)$. The series \eqref{smuuyadual} can be regarded
as a Yangian extension of the quantum immanants of \cite{o:qi}.
In particular, by
the argument of \cite[Sec.~3.4]{o:qi},
this series
does not depend on the standard tableau $\Uc$ of shape $\mu$
thus justifying the notation.

\bth\label{thm:tmuinva}
All coefficients of the series $\TT^+_{\mu}(u)\tss\vac$ belong to the
subspace of invariants $\z(\wh\Vc_{\text{cri}})$ of the extended vacuum module.
\eth

\bpf
Consider the
tensor product algebra
\beql{tenpvacy}
\underbrace{\End\CC^{N}\ot\dots\ot\End\CC^{N}}_{m+1}{}\ot\wh\Y^+(\gl_N)
\eeq
with the $m+1$ copies of $\End\CC^N$ labeled by $0,1,\dots,m$.
We need to verify the identity
\beql{aanihy}
T_0(z)\ts
\TT^+_{\mu}(u)\tss\vac=\TT^+_{\mu}(u)\tss\vac,
\eeq
where we identify the vector spaces $\wh\Vc_{\text{cri}}\cong \wh\Y^+(\gl_N)$.
By the defining relations \eqref{RTTdou}, for all $a=1,\dots,m$ we can write
\ben
T_0(z)\ts T^+_a(u+c_a)=\overline R_{0\tss a}(z-u-c_a-N/2)^{-1}\ts T^+_a(u+c_a)\ts T_0(z)
\ts \overline R_{0\tss a}(z-u-c_a+N/2).
\een
Hence, suppressing the arguments of the $R$-matrices we get
\ben
\bal
T_0(z)\ts
\tr^{}_{1,\dots,m}\ts &\Ec^{}_{\Uc}\ts T^+_1(u+c_1)\dots T^+_m(u+c_m)\tss\vac\\[0.5em]
{}&=\tr^{}_{1,\dots,m}\ts \Ec^{}_{\Uc}\ts\overline R_{0\tss 1}^{\ts-1}\dots\overline R_{0\tss m}^{\ts-1}
\ts T^+_1(u+c_1)\dots T^+_m(u+c_m)\ts T_0(z)\ts
\overline R_{0\tss m}\dots \overline R_{0\tss 1}\tss\vac\\[0.5em]
{}&=\tr^{}_{1,\dots,m}\ts \Ec^{}_{\Uc}\ts\overline R_{0\tss 1}^{\ts-1}\dots\overline R_{0\tss m}^{\ts-1}
\ts T^+_1(u+c_1)\dots T^+_m(u+c_m)\ts
\overline R_{0\tss m}\dots \overline R_{0\tss 1}\tss\vac,
\eal
\een
where the last equality holds since $T_0(z)$ acts as the identity operator on the subspace
$\End(\CC^N)^{\ot\tss (m+1)}\ot\vac$.
Using the notation \eqref{rtrain}
we get
\begin{multline}\label{rrr}
R(u_1,\dots,u_m)\tss R_{0\tss m}(u_0-u_m)\dots R_{0\tss 1}(u_0-u_1)\\[0.5em]
{}=R_{0\tss 1}(u_0-u_1)\dots R_{0\tss m}(u_0-u_m)\tss R(u_1,\dots,u_m),
\end{multline}
where $u_0$ is another variable. This follows by a repeated application of the
Yang--Baxter equation
\beql{yberep}
R_{12}(u)\ts R_{13}(u+v)\ts R_{23}(v)
=R_{23}(v)\ts R_{13}(u+v)\ts R_{12}(u)
\eeq
satisfied by the Yang $R$-matrix.
Relation \eqref{rrr} will remain valid
if each factor $R_{0\tss a}(u_0-u_a)$ is replaced with $\overline R_{0\tss a}(u_0-u_a)$.
Hence, by the fusion procedure \eqref{fusion},
the consecutive evaluations $u_a=c_a$ for $a=1,\dots,m$ imply
\ben
\Ec^{}_{\Uc}\ts\overline R_{0\tss m}(u_0-c_m)\dots \overline R_{0\tss 1}(u_0-c_1)
=\overline R_{0\tss 1}(u_0-c_1)\dots\overline R_{0\tss m}(u_0-c_m)\ts\Ec^{}_{\Uc}.
\een
By inverting the $R$-matrices we also get
\ben
\Ec^{}_{\Uc}\ts\overline R_{0\tss 1}(u_0-c_1)^{\ts-1}
\dots\overline R_{0\tss m}(u_0-c_m)^{\ts-1}
=\overline R_{0\tss m}(u_0-c_m)^{\ts-1}\dots
\overline R_{0\tss 1}(u_0-c_1)^{\ts-1}\ts\Ec^{}_{\Uc}.
\een

Returning now to the calculation of $T_0(z)\ts
\TT^+_{\mu}(u)\tss\vac$, recall that $\Ec^{}_{\Uc}$ is an idempotent and use the cyclic
property of trace together with \eqref{rtraintt} and \eqref{fusion} to write
\begin{multline}
\tr^{}_{1,\dots,m}\tss \Ec^{}_{\Uc}\tss X\tss Y
=\tr^{}_{1,\dots,m}\tss X^{\circ}\tss\Ec^{}_{\Uc}\tss Y
=\tr^{}_{1,\dots,m}\tss X^{\circ}\tss\Ec^{2}_{\Uc}\tss Y\\[0.5em]
{}=\tr^{}_{1,\dots,m}\tss \Ec^{}_{\Uc}\tss X\tss Y^{\circ}\tss\Ec^{}_{\Uc}
=\tr^{}_{1,\dots,m}\tss X\tss Y^{\circ}\tss\Ec^{2}_{\Uc}
=\tr^{}_{1,\dots,m}\tss X\tss Y^{\circ}\tss\Ec^{}_{\Uc}
=\tr^{}_{1,\dots,m}\tss X\tss\Ec^{}_{\Uc}\tss Y,
\non
\end{multline}
where we set
\ben
X=\overline R_{0\tss 1}^{\ts-1}\dots\overline R_{0\tss m}^{\ts-1},\qquad
Y=T^+_1(u+c_1)\dots T^+_m(u+c_m)\ts
\overline R_{0\tss m}\dots \overline R_{0\tss 1}
\een
and used the notation $X^{\circ}$ and $Y^{\circ}$ for the same products
written in the opposite order.
Thus, we can write
\ben
T_0(z)\ts\TT^+_{\mu}(u)\tss\vac=\tr^{}_{1,\dots,m}\tss X\tss\Ec^{}_{\Uc}\tss Y\tss\vac
=\tr^{}_{1,\dots,m}\tss X^{t_1\dots t_m}\tss\big(\Ec^{}_{\Uc}\tss Y\big)^{t_1\dots t_m}\tss\vac.
\een
We have
\ben
\big(\Ec^{}_{\Uc}\tss Y\big)^{t_1\dots t_m}=
\overline R_{0\tss m}^{\ts t_m}\dots \overline R_{0\tss 1}^{\ts t_1}
\tss\big(\Ec^{}_{\Uc}\tss T^+_1(u+c_1)\dots T^+_m(u+c_m)\big)^{t_1\dots t_m}
\een
and
\ben
X^{t_1\dots t_m}=\Big(\overline R_{0\tss 1}^{\ts-1}\Big)^{t_1}\dots
\Big(\overline R_{0\tss m}^{\ts-1}\Big)^{t_m}.
\een
By the crossing symmetry \eqref{cs}, we have
\ben
\Big(\overline R_{0\tss a}^{\ts-1}\Big)^{t_a}\ts \overline R_{0\tss a}^{\ts t_a}=1
\een
for all $a=1,\dots,m$ and so
\ben
T_0(z)\ts\TT^+_{\mu}(u)\tss\vac=\tr^{}_{1,\dots,m}\tss
\big(\Ec^{}_{\Uc}\tss T^+_1(u+c_1)\dots T^+_m(u+c_m)\big)^{t_1\dots t_m}\tss\vac
=\TT^+_{\mu}(u)\tss\vac
\een
as required.
\epf

Note two important particular cases of Theorem~\ref{thm:tmuinva} where $\mu$ is a row
or column diagram. In each case there is a unique standard tableau $\Uc$, and
the corresponding idempotents $\Ec^{}_{\Uc}$ coincide with the respective images $H^{(m)}$ and $A^{(m)}$
of the symmetrizer and anti-symmetrizer
\ben
h^{(m)}=\frac{1}{m!}\sum_{s\in\Sym_m} s\Fand
a^{(m)}=\frac{1}{m!}\sum_{s\in\Sym_m} \sgn s\cdot s
\een
under the action of $\Sym_m$ on $(\CC^N)^{\ot\ts m}$.

\bco\label{cor:symasym}
All coefficients of the series
\ben
\tr^{}_{1,\dots,m}\ts H^{(m)}\ts T^+_1(u)\dots T^+_m(u+m-1)\Fand
\tr^{}_{1,\dots,m}\ts A^{(m)}\ts T^+_1(u)\dots T^+_m(u-m+1)
\een
belong to the
subspace of invariants $\z(\wh\Vc_{\text{cri}})$ of the extended vacuum module.
\qed
\eco
In the particular case $m=N$ the second series coincides with the {\em quantum determinant}
$\qdet T^+(u)$ of the matrix $T^+(u)$; see also Proposition~\ref{prop:qdet} below.

One more family of elements of $\z(\wh\Vc_{\text{cri}})$ can be constructed
by making use of the well-known fact that the matrix
$M=T^+(u)\tss e^{-\di_u}$ with entries in the extended algebra
$\wh\Y^+(\gl_N)[[u,\di_u]]$ is a {\em Manin matrix}; see \cite{cf:mm}.
Namely, by the Newton identity \cite[Theorem~4]{cf:mm}, we have
\beql{newtonthm}
\di_z\ts \cdet(1+z\tss M)
=\cdet(1+z\tss M)\tss\sum_{m=0}^{\infty} (-z)^m\ts\tr\ts M^{m+1},
\eeq
where
\ben
\cdet(1+z\tss M)=\sum_{m=0}^{N} z^m\ts
\tr^{}_{1,\dots,m}\ts A^{(m)} M_1\dots M_m.
\een

\bco\label{cor:newton}
All coefficients of the series
\ben
\tr\ts T^+(u)\dots T^+(u-m+1),\qquad m\geqslant 1,
\een
belong to the
subspace of invariants $\z(\wh\Vc_{\text{cri}})$ of the extended vacuum module.
\eco

\bpf
Note that
\ben
M^m=T^+(u)\dots T^+(u-m+1)\ts e^{-m\tss\di_u}
\een
so that the claim follows from \eqref{newtonthm}.
\epf

\bre\label{rem:macmahon}
The MacMahon Master Theorem for Manin matrices \cite{glz:qm}
implies a relationship between the two families of Corollary~\ref{cor:symasym}:
\beql{macmahon}
\big[\cdet(1-z\tss M)\big]^{-1}=\sum_{m=0}^{\infty}
z^m\ts\tr^{}_{1,\dots,m}\ts H^{(m)} M_1\dots M_m
\eeq
for $M=T^+(u)\tss e^{-\di_u}$; see also \cite{mr:mm} for another proof
and a super-extension.
\qed
\ere

We have the following well-known
properties of quantum determinants.

\bpr\label{prop:qdet}
The coefficients of the quantum determinants
\begin{align}\label{qdet}
\qdet T(u)&=\sum_{\si\in\Sym_N}\sgn \si\cdot
t_{ \si(1)\tss 1}(u)\dots t_{ \si(N)\tss N}(u-N+1),\\
\qdet T^+(u)&=\sum_{\si\in\Sym_N}\sgn \si\cdot
t^+_{ \si(1)\tss 1}(u)\dots t^+_{ \si(N)\tss N}(u-N+1),
\label{qdetplus}
\end{align}
belong to the center of the extended double Yangian
$\DY^{\circ}(\gl_N)$.\footnote{The corresponding result for the quantum affine algebra $\U_q(\wh\gl_n)$
stated in \cite[Lemma~4.3]{fjmr:hs} holds for an arbitrary level as well (not just for the critical level).
This follows from the property $f(x)f(xq^2)\dots f(xq^{2n-2})=(1-x)/(1-xq^{2n-2})$ of the series $f(x)$
used in the proof of that lemma.}
\epr

\bpf
The respective coefficients of $\qdet T(u)$ and $\qdet T^+(u)$
are central in the Yangian $\Y(\gl_N)$ and completed dual Yangian $\wh\Y^+(\gl_N)$; see e.g.
\cite[Ch.~1]{m:yc}. Furthermore, in
the algebra \eqref{tenprkeyadual} with $m=N$ we have
\beql{rmatdet}
A^{(N)}\tss T^+_1(u)\dots T^+_N(u-N+1)=A^{(N)}\tss \qdet T^+(u).
\eeq
Arguing as in the beginning of the proof of Theorem~\ref{thm:tmuinva},
and keeping the notation, we find
\begin{multline}
T_0(z)\ts
A^{(N)}\tss T^+_1(u)\dots T^+_N(u-N+1)\\[0.5em]
{}=A^{(N)}\tss\overline R_{0\tss 1}^{\ts-1}\dots\overline R_{0\tss N}^{\ts-1}
\ts T^+_1(u)\dots T^+_N(u-N+1)\ts T_0(z)\ts
\overline R_{0\tss N}\dots \overline R_{0\tss 1}.
\non
\end{multline}
Now use the identity
\ben
A^{(N)}\tss \overline R_{0\tss N}(v+N-1)\dots \overline R_{0\tss 1}(v)=A^{(N)}.
\een
It is implied by \eqref{guprod} and the following property of the Yang $R$-matrix
\eqref{yangr}
\ben
A^{(N)}\tss R_{0\tss N}(v+N-1)\dots R_{0\tss 1}(v)=A^{(N)}\ts (1-v^{-1});
\een
see e.g. \cite[Ch.~1]{m:yc}. This proves that $T_0(z)$ commutes with $\qdet T^+(u)$.
By the same calculation, $T^+_0(z)$ commutes with $\qdet T(u)$.
\epf

Recall that the {\em vacuum module at the critical level} $V_{\text{cri}}=V_{-N}(\gl_N)$
over the affine Kac--Moody algebra $\wh\gl_N$ is defined as the quotient of $\U(\wh\gl_N)$
by the left ideal generated by $\gl_N[t]$ and $K+N$. The {\em Feigin--Frenkel center}
is the subspace
$\z(\wh\gl_N)$ of invariants
\beql{ffvac}
\z(\wh\gl_N)=\{v\in V_{\text{cri}}\ |\ \gl_N[t]\tss v=0\}.
\eeq
This subspace is a commutative associative algebra
which can be identified
with a subalgebra of $\U\big(t^{-1}\gl_N[t^{-1}]\big)$.
By a theorem of Feigin and Frenkel~\cite{ff:ak},
$\z(\wh\gl_N)$ is an algebra of polynomials in infinitely many variables; see \cite{f:lc} for
a detailed exposition of these results.

Our goal now is to use a classical limit to reproduce a construction
of elements of $\z(\wh\gl_N)$; cf.~\cite{t:qg}. By Theorem~\ref{thm:pbw},
we can regard elements of the completed dual Yangian $\wh\Y^+(\gl_N)$
as infinite linear combinations
\ben
\sum A_{i_1j_1\dots\ts i_pj_p}^{r_1\dots\ts r_p}\tss t^{(r_1)}_{i_1j_1}\dots t^{(r_p)}_{i_pj_p}
\een
of ordered monomials
over $\CC$, where all $r_i$ are negative integers.
This shows that the isomorphism \eqref{grisom} will hold true in the same form,
if we replace $\Y^+(\gl_N)$ with the completed dual Yangian $\wh\Y^+(\gl_N)$
equipped with the inherited ascending filtration defined by $\deg t_{ij}^{(-r)}=-r$; see
Corollary~\ref{cor:isompbw}. Moreover, for any element
$S\in\z(\wh\Vc_{\text{cri}})$ its image $\overline S$ in the graded algebra,
regarded as an element of $\U\big(t^{-1}\gl_N[t^{-1}]\big)$,
belongs to the Feigin--Frenkel center $\z(\wh\gl_N)$.
We will use Corollary~\ref{cor:symasym} to construct appropriate linear combinations
of elements of $\z(\wh\Vc_{\text{cri}})$ whose graded images
will be generators of $\z(\wh\gl_N)$.

Extend the ascending filtration on the completed dual Yangian to
the algebra of formal series $\wh\Y^+(\gl_N)[[u,\di_u]]$ by setting
$\deg u=1$ and $\deg\di_u=-1$ so that the associated graded algebra is
isomorphic to $\U\big(t^{-1}\gl_N[t^{-1}]\big)[[u,\di_u]]$.
Then the element
\beql{onemtd}
\tr^{}_{1,\dots,m}\tss A^{(m)} \big(1-T^+_1(u)\tss e^{-\di_u}\big)\dots
\big(1-T^+_m(u)\tss e^{-\di_u}\big)
\eeq
has degree $-m$ and its image in the graded algebra coincides with
\beql{gradep}
\tr^{}_{1,\dots,m}\tss A^{(m)} \big(\di_u+E_+(u)_1\big)\dots
\big(\di_u+E_+(u)_m\big),
\eeq
where
\beql{epu}
E_+(u)=\sum_{r=1}^{\infty} E[-r]\tss u^{r-1}.
\eeq
On the other hand, the element \eqref{onemtd} equals
\ben
\tr^{}_{1,\dots,m}\ts A^{(m)}\ts \sum_{k=0}^m\sum_{1\leqslant i_1<\dots<i_k\leqslant m}
(-1)^k\ts T^+_{i_1}(u)\dots T^+_{i_k}(u-k+1)\tss e^{-k\tss\di_u}.
\een
Transform this expression by
applying conjugations by suitable elements of $\Sym_m$ and
using the cyclic property of trace to bring it to the form
\ben
\tr^{}_{1,\dots,m}\ts A^{(m)}\ts \sum_{k=0}^m (-1)^k\binom{m}{k}\ts
T^+_{1}(u)\dots T^+_{k}(u-k+1)\tss e^{-k\tss\di_u}.
\een
Calculating partial traces of the anti-symmetrizer,
we can write this as
\beql{parteq}
\sum_{k=0}^m(-1)^k \binom{N-k}{m-k}
\tr^{}_{1,\dots,k}\ts A^{(k)}\ts
T^+_{1}(u)\dots T^+_{k}(u-k+1)\tss e^{-k\tss\di_u}.
\eeq
By Corollary~\ref{cor:symasym}, we can conclude that
all coefficients of
\eqref{gradep} belong to $\z(\wh\gl_N)$.

Together with a similar argument for the other two families
of invariants in Corollaries~\ref{cor:symasym} and \ref{cor:newton},
we thus reproduce the following result on generators of $\z(\wh\gl_N)$
from \cite{cm:ho}, \cite{ct:qs} and \cite{mr:mm};
see Corollary~\ref{cor:ffah} below.
Alternatively, for those two families it can also be derived with the use of the observation
that $M=\di_u+E_+(u)$ is a Manin matrix and applying \eqref{newtonthm}
and \eqref{macmahon}.
Introduce the power series
$\phi^{}_{m\tss a}(u)$, $\psi^{}_{m\tss a}(u)$ and $\theta^{}_{m\tss a}(u)$
by the expansions:
\begin{align}\label{deftraz}
\tr^{}_{1,\dots,m}\tss A^{(m)} \big(\di_u+E_+(u)_1\big)\dots
\big(\di_u+E_+(u)_m\big)
&=\phi^{}_{m\tss0}(u)\ts\di_u^{\tss m}
+\dots+\phi^{}_{m\tss m}(u),\\[0.7em]
\label{deftrhz}
\tr^{}_{1,\dots,m}\tss H^{(m)} \big(\di_u+E_+(u)_1\big)\dots
\big(\di_u+E_+(u)_m\big)
&=\psi^{}_{m\tss0}(u)\ts\di_u^{\tss m}
+\dots+\psi^{}_{m\tss m}(u),
\end{align}
and
\beql{deftracepa}
\tr\ts \big(\di_u+E_+(u)\big)^m=\theta^{}_{m\tss0}(u)\ts\di_u^m+\theta^{}_{m\tss1}(u)\ts\di_u^{m-1}
+\dots+\theta^{}_{m\tss m}(u).
\eeq
Define their coefficients by
\ben
\phi^{}_{m\tss m}(u)=\sum_{r=0}^{\infty}\ts \phi^{(r)}_{m\tss m}\tss u^r,\qquad
\psi^{}_{m\tss m}(u)=\sum_{r=0}^{\infty}\ts \psi^{(r)}_{m\tss m}\tss u^r
\Fand
\theta^{}_{m\tss m}(u)=\sum_{r=0}^{\infty}\ts \theta^{(r)}_{m\tss m}\tss u^r.
\een

\bco\label{cor:ffah}
Each family $\phi^{(r)}_{m\tss m}$, $\psi^{(r)}_{m\tss m}$
and $\theta^{(r)}_{m\tss m}$
with $m=1,\dots,N$ and $r=0,1,\dots$ is algebraically independent
and generates the algebra $\z(\wh\gl_N)$.
\eco

\bpf
The algebraic independence follows by considering the symbols of the elements
in the symmetric algebra as in \cite[Ch.~3]{f:lc}.
\epf

\section{Quantum vertex algebras}
\label{sec:qva}
\setcounter{equation}{0}

We will follow \cite{ek:ql5} to introduce quantum vertex algebras. We will be most concerned
with the center of a quantum vertex algebra which we introduce by analogy with vertex
algebra theory.
Our goal is to use the constructions
of invariants of the extended vacuum module given in Sec.~\ref{subsec:ivmdy} to
describe the structure of the center; see Sec.~\ref{sec:qava} below.

\subsection{Definition and basic properties}
Let $h$ be a formal parameter, $V_0$ a complex vector space and  $V=V_0 [[h]]$
a topologically free $\mathbb{C}[[h]]$-module.
Denote by $V_h ((z))$ the space of all Laurent series
\ben
v(z)=\sum_{r\in\mathbb{Z}} v_r z^{-r-1} \in V[[z^{\pm 1}]]
\een
satisfying $v_r\to 0$ as $r\to\infty$, in the $h$-adic topology.
More precisely,
$V_h ((z))$ consists of all Laurent series $v(z)$ satisfying the following condition:
for every $n\in\mathbb{Z}_{\geqslant 0}$ there exists $s\in\mathbb{Z}$ such that $r\geqslant s$ implies
$v_r \in h^n V$. Note that the space $V_h ((z))$ can be identified with $V_0 ((z))[[h]]$.

\bde\label{def:qvoa}
Let $V$ be a topologically free $\mathbb{C}[[h]]$-module.
A {\em quantum vertex algebra} $V$ over $\mathbb{C}[[h]]$ is the following data.
\begin{enumerate}[(a)]
\item\label{101} A $\mathbb{C}[[h]]$-module map (the {\em vertex operators})
\beql{truncation}
Y\,\colon\, V\otimes V \,\to\, V_h ((z)),\qquad v\ot w\mapsto Y(z)\tss (v\ot w).
\eeq
For any $v\in V$ the map $Y(v,z):V\to V_h ((z))$ is then defined by
\ben
Y(v,z)\tss w =\, Y(z)\tss (v\ot w)
\een
and which satisfies the {\em weak associativity property}:
for any $u,v,w\in V$ and $n\in\mathbb{Z}_{\geqslant 0}$
there exists $\ell\in\mathbb{Z}_{\geqslant 0}$
such that
\begin{equation}\label{associativity}
(z_0 +z_2)^\ell\ts Y(v,z_0 +z_2)Y(w,z_2)\tss u - (z_0 +z_2)^\ell\ts Y\big(Y(v,z_0)w,z_2\big)\tss u
\in h^n V[[z_0^{\pm 1},z_2^{\pm 1}]].
\end{equation}
\item\label{102} A vector $\vac\in V$ (the {\em vacuum vector}) which satisfies
\beql{v1}
Y(\vac ,z)v=v\quad\text{for all }v\in V,
\eeq
and for any $v\in V$ the series $Y(v,z)\tss\vac$ is a Taylor series in $z$ with
the property
\beql{v2}
Y(v,z)\tss\vac\big|^{}_{z=0} =v.
\eeq
\item\label{103} A $\mathbb{C}[[h]]$-module map $D\colon V\to V$
(the {\em translation operator}) which satisfies
\begin{align}
&D\tss\vac =0;\label{d1}\\
&\frac{d}{dz}Y(v,z)=[D,Y(v,z)]\quad\text{for all }v\in V.
\label{d2}
\end{align}
\item\label{104} A $\mathbb{C}[[h]]$-module map
$\mathcal{S}\colon V\otimes V\to V\otimes V\otimes\mathbb{C}((z))$ which satisfies
\begin{align}
&[D\otimes 1, \mathcal{S}(z)]=-\frac{d}{dz}\mathcal{S}(z),\label{s1}\\
\intertext{the {\em Yang--Baxter equation}}
&\mathcal{S}_{12}(z_1)\tss\mathcal{S}_{13}(z_1+z_2)\tss\mathcal{S}_{23}(z_2)
=\mathcal{S}_{23}(z_2)\tss\mathcal{S}_{13}(z_1+z_2)\tss\mathcal{S}_{12}(z_1),\label{s2}\\
\intertext{the {\em unitarity condition}}
&\mathcal{S}_{21}(z)=\mathcal{S}^{-1}(-z),\label{s3}
\end{align}
and the $\mathcal{S}$-{\em locality}:
for any $v,w\in V$ and $n\in\mathbb{Z}_{\geqslant 0}$ there exists
$\ell\in\mathbb{Z}_{\geqslant 0}$ such that for any $u\in V$
\begin{align}
&(z_1-z_2)^{\ell}\ts Y(z_1)\big(1\otimes Y(z_2)\big)\big(\mathcal{S}(z_1 -z_2)(v\otimes w)\otimes u\big)
\nonumber\\[0.4em]
&\qquad-(z_1-z_2)^{\ell}\ts Y(z_2)\big(1\otimes Y(z_1)\big)(w\otimes v\otimes u)
\in h^n V[[z_1^{\pm 1},z_2^{\pm 1}]].\label{locality}
\end{align}
\end{enumerate}
\ede

In the relations used in Definition~\ref{def:qvoa}
we applied a common expansion convention: for $\ell<0$ an expression of the form
$(z\pm w)^\ell$ should be expanded into a Taylor series of the
variable appearing on the right. For example,
\ben
(z-w)^{-1}=\sum_{r\geqslant 0}\frac{w^r}{z^{r+1}}\in\mathbb{C}((z))[[w]]
\Fand
(w-z)^{-1}=\sum_{r\geqslant 0}\frac{z^r}{w^{r+1}}\in\mathbb{C}((w))[[z]].
\een
We will apply this convention
throughout the paper, unless stated otherwise.
Also, the
tensor products are understood as
$h$-adically completed. In particular, $V\otimes V$ denotes the space
$(V_0 \otimes V_0)[[h]]$
and $V\otimes V\otimes\mathbb{C}((z))$ denotes the space
$\big(V_0\otimes V_0\otimes\mathbb{C}((z))\big)[[h]]$.

For any $r\in\ZZ$ the $r$-{\em product} $v_r w$ of elements $v$ and $w$ of $V$
is defined as the Laurent coefficient of the series
\ben
Y(v,z)\tss w =\, Y(z)\tss (v\ot w)=\sum_{r\in\mathbb{Z}}(v_r w)\ts z^{-r-1}.
\een

\bre
In the original definition of the quantum VOA in \cite{ek:ql5},
the {\em hexagon identity}
\begin{equation}\label{hexagon}
\mathcal{S}(z_1)\big(Y(z_2)\otimes 1\big)=
\big(Y(z_2)\otimes 1\big)\mathcal{S}_{23}(z_1)\mathcal{S}_{13}(z_2 +z_1)
\end{equation}
was considered instead of
the weak associativity property \eqref{associativity}.
It was proved therein that in a quantum VOA \eqref{hexagon} implies \eqref{associativity}.
Furthermore, the authors introduced the notion of {\em nondegenerate vertex algebra}
and proved that if the other axioms of quantum VOA hold, then
the hexagon identity is equivalent to the weak associativity property when
$V/hV$ is a nondegenerate vertex  algebra.
\ere

\bre
A similar notion of the $\hbar$-{\em adic {\em(}weak\ts{\em)} quantum vertex algebra} was studied
in \cite{l:hq}.
The author proved that weak associativity \eqref{associativity}
and  $\mathcal{S}$-locality \eqref{locality}
imply
the {\em Jacobi identity}
\begin{align}
&z_{0}^{-1}\delta\left(\frac{z_2 -z_1}{z_0}\right)Y(w,z_2)Y(v,z_1)\tss u\nonumber\\
&\quad-z_{0}^{-1}\delta\left(\frac{z_1 -z_2}{-z_0}\right)
Y(z_1)\big(1\otimes Y(z_2)\big)\big(\mathcal{S}(-z_0)(v\otimes w)\otimes u\big)\nonumber\\
&\quad\quad=z_1^{-1} \delta\left(\frac{z_2 -z_0}{z_1}\right)
Y\big(Y(w,z_0)v,z_1\big)\tss u\label{jacobi}
\end{align}
so it
holds for any elements $u,v,w$ of a quantum vertex algebra $V$. On the other hand,
the $\mathcal{S}$-locality and weak associativity can be recovered from \eqref{jacobi}
by using properties of the formal $\delta$-function defined by
\ben
\delta(z)=\sum_{r\in\mathbb{Z}}z^r.
\een
In particular, since $Y(w,z_0)v\in V_h ((z_0))$, for every $n\geqslant 0$
there exists $\ell\geqslant 0$ such that
$z_0^\ell\tss Y(w,z_0)\tss v$ is a Taylor series in $z_0$ modulo
$h^n V[[z_0^{\pm 1}]]$. By taking the residue $\Res_{z_0}z_0^\ell$
in \eqref{jacobi} we recover the $\mathcal{S}$-locality \eqref{locality}.
\ere

\bre\label{remark-new}
The original definition of the quantum VOA in \cite{ek:ql5} requires that
the operator $\mathcal{S}(z)$ satisfies the following property:
\begin{equation}\label{s0}
\mathcal{S}(z)(v\otimes w)-v\otimes w\ot 1 \,\in\, h \,V\otimes V\otimes
\mathbb{C}((z))\quad \text{for } v,w\in V.\end{equation}
Consequently, the classical limit $V/hV$ of a quantum VOA $V$ is a vertex algebra, as defined,
e.g. in \cite{fb:va}, \cite{flm:vo} and \cite{k:va}.
As with the definition of  $\hbar$-adic quantum vertex algebra in \cite{l:hq},
in Definition~\ref{def:qvoa} we also omit \eqref{s0}. Hence, in general,
the classical limit $V/hV$ of a quantum vertex algebra $V$ does not need
to be a vertex algebra. However, the classical limit $V/hV$  possesses
the structure of a {\em nonlocal vertex algebra} over $\mathbb{C}$, as defined in \cite{Li}.
\ere

As with the vertex algebra theory, the translation
operator $D$ is determined by the vertex operators.
Namely,
suppose that $v$ and $w$ are elements of a quantum vertex algebra $V$.
By applying \eqref{d2} to the vector $w$ and considering the coefficient of $z^{-r-1}$ we get
\begin{equation}\label{der:01}
-r\tss v_{r-1}w = D\tss v_rw -v_r\tss Dw \qquad\text{for all }r\in\mathbb{Z}.
\end{equation}
Now taking $w=\vac$, $r=-1$ in \eqref{der:01}  and using \eqref{v2} and \eqref{d1} we obtain
\begin{equation}\label{der:02}
D\tss v=v_{-2}\vac\quad\text{for all }v\in V.
\end{equation}

\subsection{The center of a quantum vertex algebra}

Let $\overline V$ be a vertex algebra. Recall that the {\em center} of $\overline V$ is defined by
\begin{equation}\label{centerofv}
\mathfrak{z}(\overline V)=\left\{v\in \overline V\ |\  w_r v=0 \text{\ \  for all  }
w\in \overline V\text{ and all } r\geqslant 0\right\};
\end{equation}
see, e.g., \cite{f:lc} and \cite{fb:va}. Equivalently, in terms of the vertex operators we have
\ben
\mathfrak{z}(\overline V)=\left\{v\in \overline V\ |\ [Y(v,z_1),Y(w,z_2)]=0
\text{\ \  for all }w\in \overline V\right\}.
\een
Consequently, the center of a vertex algebra has a structure of
a unital commutative associative algebra equipped with a derivation.
The multiplication is
defined by the $(-1)$-product
$v\cdot w= v^{}_{-1}w$ for  $v,w\in \overline V$.

We will now introduce a quantum version of $\mathfrak{z}(\overline V)$.
Let $V$ be a quantum vertex algebra. Define the {\em center} of $V$ as
the subspace
\begin{equation}\label{center}
\mathfrak{z}(V)=\left\{v\in V\ |\  w_r v=0 \text{\ \  for all  }
w\in V\text{ and all } r\geqslant 0\right\}.
\end{equation}

\bpr\label{associativity_strong_pro}
Let $V$ be a quantum vertex algebra. For any element $v\in V$ and for any
$w,u\in \mathfrak{z}(V)$ we have
\begin{equation}\label{associativity_strong}
 Y(v,z_0 +z_2)\tss Y(w,z_2)\tss u =  Y\big(Y(v,z_0)w,z_2\big)\tss u.
\end{equation}
\epr

\bpf
By \eqref{associativity} for every $n\geqslant 0$ there exists $\ell\geqslant 0$ such that
\begin{equation}\label{as:01}
(z_0 +z_2)^\ell\tss Y(v,z_0 +z_2)Y(w,z_2)u - (z_0 +z_2)^\ell\tss
Y\big(Y(v,z_0)w,z_2\big)\tss u \in h^n V[[z_0^{\pm 1},z_2^{\pm 1}]].
\end{equation}
Set
\begin{align*}
A(z_0,z_2)=Y(v,z_0 +z_2)Y(w,z_2)\tss u\Fand
B(z_0,z_2)=Y\big(Y(v,z_0)w,z_2\big)\tss u.
\end{align*}
The definition of the center \eqref{center}, together with the assumptions
$w,u\in \mathfrak{z}(V)$, imply $B(z_0,z_2)\in V[[z_0,z_2]]$. Similarly, we have
$Y(w,z_2)\tss u \in V[[z_2]]$ because $u\in \mathfrak{z}(V)$.
However, by \eqref{truncation} we have
$Y(a,z)b\in V_h ((z))$ for all $a,b\in V$, and hence
$A(z_0,z_2)\in V_0((z_0))[[h]][[z_2]]$.
Furthermore, observe that $(z_0 +z_2)^{\pm \ell}\in \mathbb{C}((z_0))((z_2))$.

Since $V_0$ is a vector space over $\mathbb{C}$, we may regard
$V_F =V_0((z_0))((h))((z_2))$ as a vector space over the field
$F=\mathbb{C}((z_0))((h))((z_2))$. By the above argument,
\ben
A(z_0,z_2),B(z_0,z_2)\in V_F\Fand (z_0 +z_2)^{\pm \ell}\in F.
\een
Therefore,  multiplying \eqref{as:01} by
$(z_0 +z_2)^{-\ell}\in \mathbb{C}((z_0))((z_2))\subset F$ we find
\begin{equation}\label{as:02}
Y(v,z_0 +z_2)Y(w,z_2)\tss u-Y\big(Y(v,z_0)w,z_2\big)\tss u \in h^n V[[z_0^{\pm 1},z_2^{\pm 1}]].
\end{equation}
Relation \eqref{as:02} holds for all $n\geqslant 0$ which implies
\begin{equation*}
 Y(v,z_0 +z_2)Y(w,z_2)\tss u -  Y\big(Y(v,z_0)w,z_2\big)\tss u=0,
\end{equation*}
as required.
\epf

We point out some consequences of Proposition~\ref{associativity_strong_pro}.
Observe that the right hand side of \eqref{associativity_strong}
is a Taylor series in the variables $z_0,z_2$:
\begin{equation}\label{as:03}
Y\big(Y(v,z_0)w,z_2\big)\tss u=\sum_{m,n<0} (v_m w)_n u\ts  z_0^{-m-1}z_2^{-n-1}.
\end{equation}
The left hand side of \eqref{associativity_strong} can be written as
\begin{align}
&Y(v,z_0 +z_2)Y(w,z_2)\tss u=\sum_{\substack{r\in\mathbb{Z}\\s<0}}
v_r w_s u \ts (z_0 +z_2)^{-r-1} z_2^{-s-1}\nonumber\\
&\qquad=\sum_{r,s<0} v_r w_s u\ts (z_0 +z_2)^{-r-1} z_2^{-s-1}
+\sum_{\substack{r\geqslant 0\\s<0}} v_r w_s u\ts (z_0 +z_2)^{-r-1} z_2^{-s-1}.\label{as:05}
\end{align}
Since the expressions in
\eqref{as:03} and \eqref{as:05} are equal by \eqref{associativity_strong}, we get
\begin{align}
&\sum_{\substack{r\geqslant 0\\s<0}} v_r w_s u (z_0 +z_2)^{-r-1}
z_2^{-s-1}=0\Fand\label{as:06}\\
&\sum_{r,s<0} v_r w_s u (z_0 +z_2)^{-r-1} z_2^{-s-1}
=\sum_{m,n<0} (v_m w)_n u z_0^{-m-1}z_2^{-n-1}.\label{as:07}
\end{align}

\bpr\label{asp:01}
The center of a quantum vertex algebra is closed under
all $s$-products with $s\in\ZZ$.
\epr

\bpf
Let $v\in V$ be an element of
a quantum vertex algebra $V$ and let $w,u\in \mathfrak{z}(V)$.
We can write \eqref{as:06} as
 \begin{equation}\label{as:08}
0=\sum_{\substack{r\geqslant 0\\s<0}} v_r w_s u (z_0 +z_2)^{-r-1} z_2^{-s-1}
=\sum_{\substack{r\geqslant 0\\s<0}}\sum_{\ell\geqslant 0}
\binom{-r-1}{\ell}v_r w_s u z_0^{-r-\ell-1} z_2^{\ell-s-1}.
\end{equation}
The coefficient of $z_0^{-1}$ in \eqref{as:08} is zero and so
 $v_0 w_s u=0$ for all $s\in\mathbb{Z}$ (note that $w_s u=0$ for $s\geqslant 0$).
By taking the coefficient of $z_0^{-2}$ we find
that $v_1 w_s u=0$ for all $s\in\mathbb{Z}$.
Continuing by the induction on
the power of $z_0$ in \eqref{as:08} we conclude that
$v_r w_s u=0$ for all $s\in\mathbb{Z}$ and $r\geqslant 0$.
Hence, $w_s u \in\mathfrak{z}(V)$ for all $s\in\mathbb{Z}$
so the proposition follows.
\epf

Define the product $\mathfrak{z}(V)\ot \mathfrak{z}(V)\to \mathfrak{z}(V)$
on the center of a quantum vertex algebra $V$
by setting
\begin{equation}\label{as:prod}
w\cdot u = w_{-1} u\qquad\text{for any }w,u\in \mathfrak{z}(V).
\end{equation}
By Proposition \ref{asp:01} the product is well-defined.

\bpr\label{prop:associativealgebra}
The product \eqref{as:prod} defines the structure of
a unital associative algebra on $\mathfrak{z}(V)$. Moreover,
this algebra is equipped with a derivation defined as
the restriction of the translation operator $D$.
\epr

\bpf
Let $v,w,u$ be arbitrary elements of $\mathfrak{z}(V)$.
By taking
the constant terms in \eqref{as:07} we get
$
(v\cdot w)\cdot u=v\cdot(w\cdot u).
$
Furthermore, \eqref{v1} implies $\vac\cdot v=v$, while \eqref{v2} implies
$v\cdot \vac=v$ and $\vac\in\mathfrak{z}(V)$, so $\vac$
is the identity in the associative algebra $\mathfrak{z}(V)$.

Taking $r\geqslant 0$ in \eqref{der:01}, we find that $Dw\in\mathfrak{z}(V)$
for any $w\in\mathfrak{z}(V)$ so that
the restriction of $D$ is a well-defined operator on
$\mathfrak{z}(V)$.
By setting $w=\vac$ and considering the
coefficient of $z_{0}$ in \eqref{as:07} we get
$v_{-2}\textstyle\vac_{-1}u=(v_{-2}\vac)_{-1}u$
for all $v,u\in \mathfrak{z}(V)$.
Therefore, using \eqref{der:02} and vacuum axiom \eqref{v1} we calculate
\ben
(Dv)\cdot u=  (Dv)_{-1}u
=(v_{-2}\vac)_{-1}u
=v_{-2}\textstyle\vac_{-1}u
=v_{-2}u.
\een
Using \eqref{der:01} with $r=-1$ we can write this as
$
Dv_{-1}u-v_{-1}Du
= D(v\cdot u) -v\cdot D(u)$.
Since $D\vac=0$ by \eqref{d1}, we conclude that
$D\colon\mathfrak{z}(V)\to\mathfrak{z}(V)$ is a derivation.
\epf

The final result of this section will demonstrate that the
center of a quantum vertex algebra is {\em $\mathcal{S}$-commutative},
as stated in the next proposition.
This replaces the commutativity property of the center of a vertex algebra
in the quantum case. In general, the center of a quantum vertex algebra
need not be commutative, as demonstrated by Proposition~\ref{prop:noncom} below.

\bpr\label{prop:scommut}
Let $V$ be a quantum vertex algebra. For any $w\in V$ and any $v,u\in\mathfrak{z}(V)$ we have
\beql{sloc}
Y(z_1)\big(1\otimes Y(z_2)\big)\big(\mathcal{S}(z_1 -z_2)(v\otimes w)\otimes u\big)
=Y(z_2)\big(1\otimes Y(z_1)\big)(w\otimes v\otimes u).
\eeq
\epr

\bpf
By \eqref{locality} for any   $w\in V$,
$v\in\mathfrak{z}(V)$ and $n\geqslant 0$ there exists
$\ell\geqslant 0$ such that for any $u\in \mathfrak{z}(V)$ we have
\begin{align}
&(z_1-z_2)^{\ell}\ts Y(z_1)\big(1\otimes Y(z_2)\big)\big(\mathcal{S}(z_1 -z_2)(v\otimes w)
\otimes u\big)\label{locality_temp_l}\\[0.4em]
&\qquad-(z_1-z_2)^{\ell}\ts Y(z_2)\big(1\otimes Y(z_1)\big)
(w\otimes v\otimes u)\in h^n V[[z_1^{\pm 1},z_2^{\pm 1}]].
\label{locality_temp}
\end{align}
Since $v$ and $u$ lie in the center of $V$ and the center is closed
under all $s$-products
by Proposition~\ref{asp:01}, the expression
$Y(z_1)(v\otimes u)$ occurring in \eqref{locality_temp} is a Taylor series in $z_1$
with coefficients in $\mathfrak{z}(V).$
Therefore, \begin{equation*}
Y(z_2)\big(1\otimes Y(z_1)\big)(w\otimes v\otimes u)\in V[[z_1,z_2]].
\end{equation*}
Now consider  \eqref{locality_temp_l}.
Recall that $\mathcal{S}(z)(a\otimes b)\in V\otimes V\otimes \mathbb{C}((z))$ for any $a,b\in V$.
Since $u\in\mathfrak{z}(V)$, the expression
$\big(1\otimes Y(z_2)\big)\big(\mathcal{S}(z_1 -z_2)(v\otimes w)\otimes u\big)$ lies
in $(V_0\otimes V_0)((z_1))[[h]][[z_2]]$ and can be written as
\begin{equation}\label{temp_exp_01}
\sum_{k\geqslant 0} \Big(\sum_{\text{fin}}
v^{(1)}_{k}\otimes v^{(2)}_{k}\otimes a_k(z_1 -z_2)\Big)  h^{k},
\end{equation}
where the internal sum is finite and
denotes an element of $V_0\otimes V_0\otimes \mathbb{C}((z_1))[[z_2]]$.
Applying the operator $Y(z_1)$ to \eqref{temp_exp_01} we get
\begin{equation}\label{temp_exp_02}
\sum_{k\geqslant 0} \Big(\sum_{\text{fin}} Y(v^{(1)}_{k},z_1) v^{(2)}_{k}
\otimes  a_k(z_1 -z_2)\Big) h^{k}.
\end{equation}
For every $m\geqslant 0$ the coefficient of $z_2^m$ in $a_k(z_1 -z_2)$
lies in $\mathbb{C}[z_1^{-1}]$,
so the internal finite sum
is an element of
$V_h ((z_1))[[z_2]]\equiv V_0 ((z_1))[[h]][[z_2]]$ for every $k\geqslant 0$.
Hence, we conclude that \eqref{temp_exp_02} lies in
$V_0 ((z_1))[[h]][[z_2]]$.
The proof is now completed as for Proposition \ref{associativity_strong_pro} by
multiplying the expression which occurs
in \eqref{locality_temp_l} and \eqref{locality_temp} by
$(z_1-z_2)^{- \ell}$.
\epf

\section{Quantum affine vertex algebra}
\label{sec:qava}
\setcounter{equation}{0}

Following \cite{ek:ql5} we introduce the quantum vertex
algebra associated with the double Yangian
for $\gl_N$.
In accordance with the general definitions of Sec.~\ref{sec:qva}, we will consider
this algebra as a module over $\CC[[h]]$. So we will start by restating
definitions of Sec.~\ref{subsec:dougln} in this context and then verify
the axioms for the quantum vertex
algebra on the vacuum module.

\subsection{Double Yangian over $\CC[[h]]$}

From now on we will work with algebras and modules over $\CC[[h]]$ and
keep the same notation for the objects associated with
the double Yangian $\DY(\gl_N)$ as in Sec.~\ref{subsec:dougln}.
The definitions of the algebras are readily translated into the $\CC[[h]]$-module context
by the formal re-scaling $u\mapsto u/h$
of the `spectral parameter' and generators
\ben
t_{ij}^{(r)}\mapsto h^{r-1}\ts t_{ij}^{(r)},\qquad t_{ij}^{(-r)}\mapsto h^{-r}\ts t_{ij}^{(-r)},
\qquad C\mapsto C,
\een
for $r\geqslant 1$. Conversely, the formal evaluation $h=1$ can be used to recover
some of the definitions and formulas of Sec.~\ref{subsec:dougln}.
The Yang $R$-matrix \eqref{yangr} now takes the form
\beql{yangrh}
R(u)=1-P\tss h\tss u^{-1},
\eeq
while for the normalized $R$-matrix \eqref{orm} we have
\beql{ormh}
\overline R(u)=g(u/h)\ts\big(1-P\tss h\tss u^{-1}\big).
\eeq

The double Yangian $\DY(\gl_N)$ is now defined as the associative algebra over $\CC[[h]]$
generated by the central element $C$ and elements
$t_{ij}^{(r)}$ and $t_{ij}^{(-r)}$, where $1\leqslant i,j\leqslant N$ and $r=1,2,\dots$,
subject to the defining relations
\begin{align}
R(u-v)\ts T_1(u)\ts T_2(v)&=T_2(v)\ts T_1(u)\ts R(u-v),
\label{RTTh}\\[0.5em]
R(u-v)\ts T^+_1(u)\ts T^+_2(v)&=T^+_2(v)\ts T^+_1(u)\ts R(u-v),
\label{RTTdualh}\\[0.5em]
\overline R\big(u-v+h\tss C/2\big)\ts T_1(u)\ts T^+_2(v)&=T^+_2(v)\ts T_1(u)\ts
\overline R\big(u-v-h\tss C/2\big),
\label{RTTdouh}
\end{align}
where the matrices $T(u)$ and $T^+(u)$ are given by
\beql{ttp}
T(u)=\sum_{i,j=1}^N e_{ij}\ot t_{ij}(u)\Fand
T^+(u)=\sum_{i,j=1}^N e_{ij}\ot t^+_{ij}(u)
\eeq
with
\ben
t_{ij}(u)=\de_{ij}+h\tss\sum_{r=1}^{\infty}t_{ij}^{(r)}\ts u^{-r}
\Fand
t^+_{ij}(u)=\de_{ij}-h\tss\sum_{r=1}^{\infty}t_{ij}^{(-r)}\ts u^{r-1}.
\een
The coproduct now takes the form
\begin{align*}
&\Delta:t_{ij}(u)\mapsto\sum_{k=1}^{N} t_{ik}\big(u+ h\tss C_2 /4\big)
\otimes t_{kj}\big(u- h\tss C_1 /4\big),\\
&\Delta:t_{ij}^{+ }(u)\mapsto\sum_{k=1}^{N} t_{ik}^{+ }\big(u- h\tss C_2 /4\big)
\otimes t_{kj}^{+ }\big(u+ h\tss C_1 /4\big),\\[0.3em]
&\Delta:C\mapsto C\otimes 1+1\otimes C,
\end{align*}
where the tensor products are understood as
$h$-adically completed. The antipode and counit are defined by the same formulas as for
the extended double Yangian $\DY^{\circ}(\gl_N)$; see Sec.~\ref{subsec:dougln}.

The Poincar\'e--Birkhoff--Witt theorem for the double Yangian extends
to the algebra $\DY(\gl_N)$ over $\CC[[h]]$; see Theorem~\ref{thm:pbw}.
Therefore the subalgebra of $\DY(\gl_N)$ generated by the elements
$t_{ij}^{(r)}$ with $1\leqslant i,j\leqslant N$ and $r\geqslant 1$ can be identified
with the Yangian $\Y(\gl_N)$ defined by the relations \eqref{RTTh}.
Similarly, the subalgebra generated by the elements
$t_{ij}^{(-r)}$ with $1\leqslant i,j\leqslant N$ and $r\geqslant 1$ can be identified
with the dual Yangian $\Y^+(\gl_N)$ defined by the relations \eqref{RTTdualh}.

\subsection{Vacuum module as a quantum vertex algebra}

The double Yangian at the level $c\in\CC$ is the quotient $\DY_c(\gl_N)$
of $\DY(\gl_N)$ by the ideal generated by $C-c$.
Similar to Sec.~\ref{subsec:ivmdy},
the {\em vacuum module $\Vc_c(\gl_N)$ at the level} $c$ over the double Yangian
is the $h$-adic completion of the quotient
\beql{vacmyh}
\DY_c(\gl_N)/\DY_c(\gl_N)\langle t_{ij}^{\ts(r)}\ts|\ts r\geqslant 1\rangle.
\eeq
By the Poincar\'e--Birkhoff--Witt theorem (Theorem~\ref{thm:pbw}),
we can identify \eqref{vacmyh} with the dual Yangian
$\Y^+(\gl_N)$ as a $\CC[[h]]$-module.

As demonstrated in \cite{ek:ql5}, the $h$-adic completion of \eqref{vacmyh}
 possesses a quantum vertex algebra
structure. In the classical limit $h\to 0$ it turns into the affine vertex algebra
$V_c(\gl_N)$.
Accordingly, $\Vc_c(\gl_N)$ is called the {\em quantum affine vertex algebra}.
To introduce the structure, we need some notation.
For a positive integer $n$,
consider the tensor product space
\beql{ev}
(\End\mathbb{C}^{N})^{\otimes n} \otimes \Vc_c(\gl_N).
\eeq
Given a variable $z$ and a
family of variables
$u=(u_1,\dots,u_n)$, set
\begin{align}
T_n(u)&=T_{1\,n+1}(u_1)\dots T_{n\,n+1}(u_n),\nonumber\\[0.5em]
T_n^{+}(u)&=T_{1\,n+1}^{+}(u_1)\dots T_{n\,n+1}^{+}(u_n),\nonumber\\[0.5em]
T_n(u|z)&=T_{1\,n+1}(z+u_1)\dots T_{n\,n+1}(z+u_n),\non\\[0.5em]
T_n^{+}(u|z)&=T_{1\,n+1}^{+}(z+u_1)\dots T_{n\,n+1}^{+}(z+u_n).\nonumber
\end{align}
Here we extend the notation \eqref{cab} to include
the vacuum module as a component of tensor products
so that the subscript $n+1$ corresponds to
$\Vc_c(\gl_N)$. For $n=0$
these products will be considered as being equal to the identity.
The respective components of the matrices \eqref{ttp}
are understood as operators on $\Vc_c(\gl_N)$.
The series $T_{i\,n+1}(z+u_i)$ and $T^+_{i\,n+1}(z+u_i)$
should be expanded in nonpositive and nonnegative powers
of $z$, respectively.

For nonnegative integers $m$ and $n$ introduce functions depending on a variable $z$ and the
families of variables
$u=(u_1,\dots,u_n)$ and $v=(v_1,\dots,v_m)$ with values in
the space
\beql{eevv}
(\End\mathbb{C}^{N})^{\otimes n} \otimes
(\End\mathbb{C}^{N})^{\otimes m}
\eeq
by
\beql{Rnm12}
R_{nm}^{12}(u|v|z)= \prod_{j=1,\dots,n}^{\longrightarrow}\ts\ts
\prod_{i=n+1,\dots,n+m}^{\longleftarrow} R_{ji}(z+u_j -v_{i-n})
\eeq
with the arrows indicating the order of the factors, where we
use the Yang $R$-matrix \eqref{yangrh} and adopt the matrix notation
as in \eqref{cab}. As above, empty products will be understood as being equal
to the identity.
We also define $\overline{R}_{nm}^{\ts 12}(u|v|z)$ by the same formula \eqref{Rnm12},
where the $R$-matrix \eqref{ormh} is used instead of $R(u)$.
The superscripts $1$ and $2$ are meant to indicate
the tensor factors in \eqref{eevv}. We also adopt the superscript notation
for multiple tensor products of the form
\beql{eeevvv}
(\End\mathbb{C}^{N})^{\otimes n} \otimes
(\End\mathbb{C}^{N})^{\otimes m}\otimes
(\End\mathbb{C}^{N})^{\otimes k}\otimes \Vc_c(\gl_N)
\otimes \Vc_c(\gl_N)\otimes \Vc_c(\gl_N).
\eeq
Expressions like $T_{n}^{14}(u)$ or $T_{k}^{ 35}(u)$ will be understood
as the respective operators $T_n(u)$ or $T_k(u)$, whose non-identity components
belong to the corresponding tensor factors. In particular, the non-identity
components of $T_{k}^{ 35}(u)$ belong to the
factors
\ben
n+m+1,\,n+m+2,\,\dots,\,n+m+k\fand n+m+k+2.
\een

Employing this notation, we point out some immediate consequences of the defining relations
\eqref{RTTh}--\eqref{RTTdouh} for operators on
\ben
(\End\mathbb{C}^{N})^{\otimes n} \otimes
(\End\mathbb{C}^{N})^{\otimes m}\otimes \Vc_c(\gl_N).
\een
They follow by a straightforward induction and take the form
\begin{align}
R_{nm}^{12}(u|v|z-w)T_n^{+13}(u|z)T_m^{+23}(v|w)
&=T_m^{+23}(v|w)T_n^{+13}(u|z)R_{nm}^{12}(u|v|z-w),\label{rtt6}\\[0.5em]
R_{nm}^{12}(u|v|z-w)T_n^{13}(u|z)T_m^{23}(v|w)
&=T_m^{23}(v|w)T_n^{13}(u|z)R_{nm}^{12}(u|v|z-w),\label{rtt7}\\[0.5em]
\overline{R}_{nm}^{\ts 12}(u|v|z-w+h\tss c/2)T_n^{13}(u|z)T_m^{+23}(v|w)
&=T_m^{+23}(v|w)T_n^{13}(u|z)\overline{R}_{nm}^{\ts 12}(u|v|z-w-h\tss c/2).
\label{rtt8}
\end{align}

It will also be convenient to use an ordered product notation
for elements of the tensor product of two associative algebras $\Ac\ot\Bc$.
Suppose that $A_1,A_2\in\Ac$ and $B_1,B_2\in\Bc$. Let
$F=A_1 \otimes B_1$ and define the following products
\begin{alignat}{2}
{}^{ll}F(A_2\otimes B_2)&=A_1 A_2 \otimes B_1 B_2,\qquad
^{lr}F(A_2\otimes B_2)&&=A_1 A_2 \otimes  B_2 B_1,
\non\\[0.3em]
^{rl}F(A_2\otimes B_2)&= A_2 A_1 \otimes B_1 B_2, \qquad
^{rr}F(A_2\otimes B_2)&&= A_2 A_1 \otimes B_2 B_1,
\label{lr}
\end{alignat}
indicating the left and right multiplication of the components.
For  $\alpha,\beta\in\{l,r\}$ we will denote by
$({}^{\alpha\beta}F)^{-1}$ the operator $G$ such that
$({}^{\alpha\beta}G)\tss F=1$. Note that $({}^{\alpha\beta}F)^{-1}$
and ${}^{\alpha\beta}\big(F^{-1}\big)$ need not be equal.

This notation will often be applied to products of
$R$-matrices
$F=R_{nm}^{12}(u|v|z)$, where the roles of $\Ac$ and $\Bc$ will be played
by the first and second components in \eqref{eevv}.
We point out the formulas for the inverse operators associated with the $R$-matrix
\eqref{yangrh}:
\begin{align*}
&\big({}^{lr}R(u)\big)^{-1}=\big({}^{rl}R(u)\big)^{-1}
=\big(1-hNu^{-1}\big)^{-1}\big(R(-u)-hNu^{-1}\big),\\[0.4em]
&\big({}^{ll}R(u)\big)^{-1}=\big({}^{rr}R(u)\big)^{-1}
=R(u)^{-1}=\big(1-h^2 u^{-2}\big)^{-1}R(-u),
\end{align*}
which can be used to calculate the inverse operators corresponding to
$F=R_{nm}^{12}(u|v|z)$.

We will now use the general definition of quantum vertex algebra reproduced
in Sec.~\ref{sec:qva}; see Definition~\ref{def:qvoa}.
The following theorem is due to Etingof and Kazhdan~\cite{ek:ql5}.

\bth\label{EK:qva}
There exists a unique well-defined structure of quantum vertex algebra
on the vacuum module $\Vc_c(\gl_N)$ with the following data.
\begin{enumerate}[(a)]
\item\label{102h} The vacuum vector is
\beql{qva2}
\vac=1\in \Vc_c(\gl_N).
\eeq
\item\label{101h} The vertex operators are defined by
\beql{qva1}
Y\big(T_n^{+}(u)\vac,z\big)=T_n^{+}(u|z)\ts T_n (u|z+h\tss c/2)^{-1}.
\eeq
\item\label{103h} The translation operator $D$ is defined by
\beql{qva3}
e^{zD}\ts T^+ (u_1)\dots T^+ (u_n)\vac = T^+ (z+u_1)\dots T^+ (z+u_n)\vac.
\eeq
\item\label{104h} The map
$\mathcal{S}$ is defined by the relation
\begin{multline}
\mathcal{S}_{34}(z)\Big(\overline{R}_{nm}^{\ts 12}(u|v|z)^{-1}\ts T_{m}^{+24}(v)\ts
\overline{R}_{nm}^{\ts 12}(u|v|z-h\tss c)\ts T_{n}^{+13}(u)(\vac\otimes \vac) \Big)\\[0.4em]
{}=T_{n}^{+13}(u)\ts \overline{R}_{nm}^{\ts 12}(u|v|z+h\tss c)^{-1}\ts
T_{m}^{+24}(v)\ts \overline{R}_{nm}^{\ts 12}(u|v|z)(\vac\otimes \vac)\label{qva4}
\end{multline}
for operators on
\beql{eevvh}
(\End\mathbb{C}^{N})^{\otimes n} \otimes
(\End\mathbb{C}^{N})^{\otimes m}\otimes \Vc_c(\gl_N)\otimes \Vc_c(\gl_N).
\eeq
\end{enumerate}
\eth

\bpf
We add some details as compared to \cite{ek:ql5}, to take care
of the variations of the definition of the quantum vertex algebra.
Let $V=\Vc_c(\gl_N)$. We start by pointing out that $Y$ is a well-defined operator
as in \eqref{truncation}. Indeed,
since the coefficients of the series $T_n^{+}(u)\vac$ span the $h$-adically dense subset of
$\Vc_c(\gl_N)$, it suffices
to verify that $Y$ preserves the ideal of relations of the dual Yangian.
This follows by employing \eqref{RTTh} and \eqref{RTTdualh} as in the proof of
\cite[Lemma~2.1]{ek:ql5}.
As a next step, we will verify
the weak associativity property \eqref{associativity}.
Let $m$, $n$ and $k$ be nonnegative integers and
let $u=(u_1,\dots,u_n)$, $v=(v_1,\dots,v_m)$ and $w=(w_1,\dots,w_k)$ be families of variables.
Note the following relation which is a consequence of \eqref{rtt8}:
\begin{multline}
T^{14}_{n}(u|z_0 +h\tss c/2)^{-1}\ts \overline{R}_{nm}^{\ts 12}(u|v|z_0+h\tss c)^{-1}
\ts T^{+24}_{m}(v)\\[0.4em]
{}=T^{+24}_{m}(v)\ts
\overline{R}_{nm}^{\ts 12}(u|v|z_0)^{-1}\ts T^{14}_{n}(u|z_0 +h\tss c/2)^{-1}.
\label{rtt9}
\end{multline}
Here and below we use the additional variables $z_0$ and $z_2$ as in \eqref{associativity}.
Using  the definition of the vacuum module together
with \eqref{qva1} and \eqref{rtt9} we get
\begin{align}
Y\big(T^{+14}_{n}(u)\vac,z_0\big)&\ts
\overline{R}_{nm}^{\ts 12}(u|v|z_0+h\tss c)^{-1}\ts T^{+24}_{m}(v)\vac
\non\\[0.4em]
&= T^{+14}_{n}(u|z_0)\ts
T^{14}_{n}(u|z_0+h\tss c/2)^{-1}\ts
\overline{R}_{nm}^{\ts 12}(u|v|z_0+h\tss c)^{-1}\ts
T^{+24}_{m}(v)\vac
\non\\[0.4em]
&=T^{+14}_{n}(u|z_0)\ts
T^{+24}_{m}(v)\ts\overline{R}_{nm}^{\ts 12}(u|v|z_0)^{-1}\ts T^{14}_{n}(u|z_0+h\tss c/2)^{-1}
\vac\non\\[0.4em]
&=T^{+14}_{n}(u|z_0)\ts
T^{+24}_{m}(v) \ts\overline{R}_{nm}^{\ts 12}(u|v|z_0)^{-1}\vac.
\non
\end{align}
For fixed positive integers $M$ and $p$
and operators $A$ and $B$ on \eqref{eeevvv} of this form,
we will say that $A$ and $B$ are {\em equivalent},
if the coefficients of all monomials
\begin{equation}\label{monomials:01}
u_1^{r_1}\dots u_{n}^{r_n}v_{1}^{s_1}\dots v_{m}^{s_m}w_{1}^{t_1}\dots w_{k}^{t_k}\qquad
\text{with}\quad 0\leqslant r_1,\dots,r_n,s_1,\dots,s_m,t_1,\dots,t_k\leqslant M
\end{equation}
in $A-B$ belong to the subspace
$h^p V[[z_0^{\pm 1},z_{2}^{\pm 1}]]$.
Let $\ell$ be a nonnegative integer such that
the coefficients of the monomials \eqref{monomials:01}
in the operator
\ben
z^\ell\ts T^{14}_{n}(u|z+h\tss c/2)^{-1}\tss T^{+34}_{k}(w)\vac
\een
have only nonnegative powers of $z$ modulo $h^p$. By
the above calculation, the operator
\beql{aas1}
(z_0 +z_2)^\ell\ts Y\Big(Y\big(T^{+14}_{n}(u)\vac,z_0\big)\ts
\overline{R}_{nm}^{\ts 12}(u|v|z_0+h\tss c)^{-1}\ts T^{+24}_{m}(v)\vac,z_2\Big)\ts
T^{+34}_{k}(w)\vac
\eeq
equals
\ben
(z_0 +z_2)^\ell\ts  Y\Big(T^{+14}_{n}(u|z_0)\ts
T^{+24}_{m}(v)\ts \overline{R}_{nm}^{\ts 12}(u|v|z_0)^{-1}\vac,z_2\Big)\ts T^{+34}_{k}(w)\vac
\een
which by \eqref{qva1} coincides with
\begin{multline}
(z_0 +z_2)^\ell\ts  T^{+14}_{n}(u|z_2+z_0)\ts
T^{+24}_{m}(v|z_2)\ts T^{24}_{m}(v|z_2+h\tss c/2)^{-1}\\[0.4em]
{}\times T^{14}_{n}(u|z_2+z_0+h\tss c/2)^{-1}\ts
\overline{R}_{nm}^{12}(u|v|z_0)^{-1}\ts T^{+34}_{k}(w)\vac.
\non
\end{multline}
By our assumption on $\ell$, only nonnegative powers of $z_0 +z_2$ will occur
in the expansion of this operator modulo $h^p$, so that we may swap $z_0$ and $z_2$
to get an equivalent operator
\begin{multline}
(z_0 +z_2)^\ell\ts  T^{+14}_{n}(u|z_0+z_2)\ts
T^{+24}_{m}(v|z_2)\ts T^{24}_{m}(v|z_2+h\tss c/2)^{-1}\\[0.4em]
{}\times T^{14}_{n}(u|z_0+z_2+h\tss c/2)^{-1}\ts
\overline{R}_{nm}^{12}(u|v|z_0)^{-1}\ts T^{+34}_{k}(w)\vac.
\label{aas2}
\end{multline}
On the other hand, by \eqref{qva1} the operator
\beql{aas3}
Y\big(T^{+14}_{n}(u)\vac,z_0 +z_2\big)\ts
Y\big(\overline{R}_{nm}^{\ts 12}(u|v|z_0 +h\tss c)^{-1} T^{+24}_{m}(v)\vac,z_2 \big)\ts T^{+34}_{k}(w)\vac
\eeq
equals
\begin{multline}
T^{+14}_{n}(u|z_0+z_2)\ts
T^{14}_{n}(u|z_0+z_2+h\tss c/2)^{-1}\\[0.4em]
{}\times\overline{R}_{nm}^{\ts 12}(u|v|z_0 +h\tss c)^{-1}\ts T^{+24}_{m}(v|z_2)\ts
T^{24}_{m}(v|z_2+h\tss c/2)^{-1}\ts T^{+34}_{k}(w)\vac.
\non
\end{multline}
Applying \eqref{rtt9} and then \eqref{rtt7} we can write this as
\begin{multline}
T^{+14}_{n}(u|z_0+z_2)\ts T^{+24}_{m}(v|z_2)\\[0.4em]
{}\times T^{24}_{m}(v|z_2+h\tss c/2)^{-1}\ts T^{14}_{n}(u|z_0+z_2+h\tss c/2)^{-1}
\ts\overline{R}_{nm}^{\ts 12}(u|v|z_0 )^{-1}\ts T^{+34}_{k}(w)\vac.
\non
\end{multline}
Observe that after multiplication by $(z_0 +z_2)^\ell$
this coincides with \eqref{aas2}. Therefore, when the operator
\eqref{aas3} is multiplied by $(z_0 +z_2)^\ell$, it will
be equal to \eqref{aas1}
modulo $h^p$.
By applying ${}^{rl}\big(({}^{rl}\overline{R}_{nm}^{\ts 12}(u|v|z_0 +h\tss c)^{-1})^{-1}\big)$
to both sides of this equality we get \eqref{associativity}, as required.

The vacuum axioms \eqref{v1} and \eqref{v2}
are immediate from the definitions of the vacuum vector
and vertex operators.

Now we verify the translation operator $D$ is well-defined
by \eqref{qva3} and satisfies
the axioms \eqref{d1} and \eqref{d2}. We need to check that $D$
preserves the defining relations \eqref{RTTdualh} of the dual Yangian.
This is a straightforward calculation; cf. \cite[Lemma~2.1]{ek:ql5}.
Furthermore, $e^{zD}\vac =\vac$ so that
\eqref{d1} holds. Now suppose that $m$ and $n$ are nonnegative integers.
Taking the coefficient of $z$ in \eqref{qva3} we get
\begin{equation}\label{eq:002}
D\ts T^+ (u_1)\dots T^+ (u_n)\vac =\Big( \sum_{l=1}^{n}
\frac{\partial}{\partial u_l}\tss\Big)\ts T^+ (u_1)\dots  T^+ (u_n)\vac.
\end{equation}
Therefore, using \eqref{qva1} we obtain
\ben
\frac{\partial}{\partial z}\ts Y\big(T_n^{+13}(u)\vac,z\big)\ts  T_{m}^{+23}(v)\vac
=\frac{\partial}{\partial z}\ts  T_n^{+13}(u|z)\ts T_n^{13}(u|z+h\tss c/2)^{-1}\ts
T_{m}^{+23}(v)\vac
\een
which can be written as
\ben
\Big(\sum_{l=1}^{n}\frac{\partial}{\partial u_l}\tss\Big)\ts
T_n^{+13}(u|z)\ts T_n^{13}(u|z+h\tss c/2)^{-1}\ts
T_{m}^{+23}(v)\vac.
\een
This coincides with
\ben
D\ts Y\big(T_n^{+13}(u)\vac,z\big)\ts T_{m}^{+23}(v)\vac
-Y\big(T_n^{+13}(u)\vac,z\big)\ts  D\ts T_{m}^{+23}(v)\vac,
\een
since
\begin{multline}
D\ts Y\big(T_n^{+13}(u)\vac,z\big)\ts T_{m}^{+23}(v)\vac\\
{}=\Big(\sum_{l=1}^{n}\frac{\partial}{\partial u_l}
+ \sum_{k=1}^{m}\frac{\partial}{\partial v_k}\tss\Big)\ts
T_n^{+13}(u|z)\ts T_n^{13}(u|z+h\tss c/2)^{-1}\ts
T_{m}^{+23}(v)\vac
\non
\end{multline}
and
\ben
Y\big(T_n^{+13}(u)\vac,z\big)\ts  D\ts T_{m}^{+23}(v)\vac
=
T_n^{+13}(u|z)\ts T_n^{13}(u|z+h\tss c/2)^{-1}
\ts\Big(\sum_{k=1}^{m}\frac{\partial}{\partial v_k}\tss\Big)\ts
T_{m}^{+23}(v)\vac,
\een
thus verifying \eqref{d2}.
Now turn to the axioms concerning the map $\Sc$.
Using the notation \eqref{lr}, we can write the operators appearing in
\eqref{qva4} in the form
\begin{multline}
\overline{R}_{nm}^{\ts 12}(u|v|z)^{-1}\ts T_{m}^{+24}(v)\ts
\overline{R}_{nm}^{\ts 12}(u|v|z-hc)\ts T_{n}^{+13}(u)\\[0.5em]
{}=
{}^{ll}\overline{R}_{nm}^{\ts 12}(u|v|z)^{-1}\ts  {}^{lr}
\overline{R}_{nm}^{\ts 12}(u|v|z-hc)\ts T_{n}^{+13}(u)T_{m}^{+24}(v)
\non
\end{multline}
and
\begin{multline}
T_{n}^{+13}(u)\ts \overline{R}_{nm}^{12}(u|v|z+hc)^{-1}\ts T_{m}^{+24}(v)
\ts \overline{R}_{nm}^{12}(u|v|z)\\[0.5em]
{}=
{}^{rr}\overline{R}_{nm}^{\ts 12}(u|v|z)\ts
{}^{rl}\overline{R}_{nm}^{\ts 12}(u|v|z+hc)^{-1}\ts
T_{n}^{+13}(u)\ts T_{m}^{+24}(v).
\non
\end{multline}
Hence \eqref{qva4} can be  written as
\begin{multline}
\mathcal{S}_{34}(z)\Big({}^{ll}\overline{R}_{nm}^{\ts 12}(u|v|z)^{-1}\ts  {}^{lr}
\overline{R}_{nm}^{\ts 12}(u|v|z-hc)
\ts T_{n}^{+13}(u)T_{m}^{+24}(v)(\vac\otimes \vac) \Big)\\[0.4em]
\quad ={}^{rr}\overline{R}_{nm}^{\ts 12}(u|v|z)\ts
{}^{rl}\overline{R}_{nm}^{\ts 12}(u|v|z+hc)^{-1}\ts
T_{n}^{+13}(u)\ts T_{m}^{+24}(v)(\vac\otimes \vac)
\non
\end{multline}
which is equivalent to
\begin{multline}
\mathcal{S}_{34}(z)\Big(T_{n}^{+13}(u)\ts T_{m}^{+24}(v)(\vac\otimes \vac)\Big)=
{}^{lr}\big(({}^{rl}\overline{R}_{nm}^{12}(u|v|z-h\tss c))^{-1}\big)\ts\ts
{}^{ll}\overline{R}_{nm}^{\ts 12}(u|v|z)\\[0.5em]
{}\times{}^{rr}\overline{R}_{nm}^{\ts 12}(u|v|z)\ts\ts{}^{rl}\overline{R}_{nm}^{\ts 12}(u|v|z+h\tss c)^{-1}
\ts T_{n}^{+13}(u)\ts T_{m}^{+24}(v)(\vac\otimes \vac).
\non
\end{multline}
This form of $\Sc$ is convenient for checking
that the map is well-defined; cf. \cite[Lemma~2.1]{ek:ql5}.
 Property
\eqref{s1} is checked in the same way as \eqref{d2} with the use of \eqref{eq:002}.
The Yang--Baxter equation \eqref{s2},
the unitarity condition \eqref{s3} and the $\Sc$-locality property \eqref{locality}
are verified by straightforward calculations
which rely on the Yang--Baxter equation \eqref{yberep} satisfied by the $R$-matrix \eqref{ormh}
and the unitarity property \eqref{up}.
\epf

We now give an example based on the structure of the dual Yangian to demonstrate that
the center of a quantum
vertex algebra need not be commutative, in general.
We use the same notation for products of generators matrices
as in the beginning of this section.

\bpr\label{prop:noncom}
There exists a unique well-defined structure of quantum vertex algebra
on the $h$-adic completion $V$ of the $\CC[[h]]$-module $\Y^+(\gl_N)$ with the following data.
\begin{enumerate}[(a)]
\item\label{102he} The vacuum vector is
\beql{qqva2}
\vac=1\in \Y^+(\gl_N).
\eeq
\item\label{101he} The vertex operators are defined by
\beql{qqva1}
Y\big(T_n^{+}(u)\vac,z\big)=T_n^{+}(u|z).
\eeq
\item\label{103he} The translation operator $D$ is defined by
\beql{qqva3}
e^{zD}\ts T^+ (u_1)\dots T^+ (u_n)\vac = T^+ (z+u_1)\dots T^+ (z+u_n)\vac.
\eeq
\item\label{104he} The map
$\mathcal{S}$ is defined by the relation
\begin{multline}
\mathcal{S}_{34}(z)\Big(T_{n}^{+13}(u)T_{m}^{+24}(v)(\vac\otimes \vac) \Big)\\[0.4em]
{}=\overline{R}_{nm}^{\ts 12}(u|v|z)\ts T_{n}^{+13}(u)\ts
T_{m}^{+24}(v)\ts \overline{R}_{nm}^{\ts 12}(u|v|z)^{-1}(\vac\otimes \vac).
\label{qqva4}
\end{multline}
\end{enumerate}

Moreover, the center $\z(V)$ of the quantum vertex algebra $V$
coincides with $V$.
\epr

\bpf
The last claim follows since
the image of the vertex operator map $Y$ is contained in
$V[[z]]$. In particular, $\z(V)$ is not commutative for $N\geqslant 2$.

The maps $Y$, $D$ and $\Sc$ are well-defined, as follows by the
same arguments as for the
proof of Theorem~\ref{EK:qva}.
The quantum vertex algebra axioms are also checked in a similar way
with some obvious modifications. We only verify
the $\Sc$-commutativity
\eqref{sloc} which implies the $\Sc$-locality property \eqref{locality}.
Set $z=z_1 -z_2$ and consider
the left hand side in \eqref{sloc}. The application of
$\mathcal{S}_{45}(z)\otimes 1$ to
\ben
T_{n}^{+14}(u)\ts  T_{m}^{+25}(v)\ts T_{k}^{+36}(w)
(\vac\otimes\vac\otimes\vac)
\een
gives
\ben
\overline{R}_{nm}^{\ts 12}(u|v|z)\ts
T_{n}^{+14}(u)\ts T_{m}^{+25}(v)\ts \overline{R}_{nm}^{\ts 12}(u|v|z)^{-1}
\ts T_{k}^{+36}(w)(\vac\otimes\vac\otimes\vac).
\een
Further applying $1\otimes Y(z_2)$ we get
\ben
\overline{R}_{nm}^{\ts 12}(u|v|z)\ts
T_{n}^{+14}(u)\ts T_{m}^{+25}(v|z_2)\ts \overline{R}_{nm}^{\ts 12}(u|v|z)^{-1}
\ts T_{k}^{+35}(w)(\vac\otimes\vac)
\een
which becomes
\beql{new:01}
\overline{R}_{nm}^{\ts 12}(u|v|z)\ts T_{n}^{+14}(u|z_1)\ts
T_{m}^{+24}(v|z_2)\ts \overline{R}_{nm}^{\ts 12}(u|v|z)^{-1}\ts T_{k}^{+34}(w)\vac
\eeq
after the application of $Y(z_1)$.
For the right hand side we have
\ben
T_{m}^{+24}(v)\ts T_{n}^{+15}(u)\ts T_{k}^{+36}(w)
(\vac\otimes\vac\otimes\vac)
\xmapsto{1\ts\otimes\ts Y(z_1)}T_{m}^{+24}(v)\ts T_{n}^{+15}(u|z_1)
\ts T_{k}^{+35}(w)(\vac\otimes\vac)
\een
and the application of $Y(z_2)$ gives
\beql{new:02}
T_{m}^{+24}(v|z_2)\ts T_{n}^{+14}(u|z_1)\ts
T_{k}^{+34}(w)\vac.
\eeq
Now \eqref{rtt6} implies that \eqref{new:01} coincides with
\eqref{new:02} and so the $\mathcal{S}$-commutativity property
\eqref{sloc} follows.
\epf

\bre
The classical limit $V/hV\equiv U(t^{-1}\mathfrak{gl}_N[t^{-1}])$ of the
quantum vertex algebra $V$ from Proposition \ref{prop:noncom}
  is a  nonlocal vertex algebra over $\mathbb{C}$; see Remark \ref{remark-new}.
Its vertex operator map is given by
\ben
Y(E_{+1}(u_1)\ldots E_{+n}(u_n)\vac,z )=E_{+1}(u_1 +z)
\dots E_{+n}(u_n +z)\quad \text{for }n\geqslant 0,
\een
where
$E_{+}(u)\in U(t^{-1}\mathfrak{gl}_N[t^{-1}])[[u]]$ is defined by \eqref{epu}.
\ere

\subsection{Central elements of the completed double Yangian}
\label{subsec:cedy}

As with the affine vertex algebras,
the vertex operator formulas \eqref{qva1} suggest a construction
of central elements of a completed double Yangian; cf. \cite[Sec.~4.3.2]{f:lc}
and Remark~\ref{rem:quace} below.
However, we will not use the quantum vertex algebra structure, but rather
give a direct proof as in \cite{fjmr:hs}.

Introduce the completion of the double Yangian $\DY_c(\gl_N)$ at the level $c$
as the inverse limit
\beql{compluady}
\wt\DY_c(\gl_N)=\lim_{\longleftarrow}\DY_c(\gl_N)/I_p,
\eeq
where $p\geqslant 1$ and $I_p$ denotes the left ideal of $\DY_c(\gl_N)$, generated
by all elements $t_{ij}^{(r)}$ with $r\geqslant p$.
Using the idempotents $\Ec^{}_{\Uc}$ as in Theorem~\ref{thm:tmuinva},
introduce the Laurent series in $u$ with
coefficients in the $h$-adically completed algebra of formal power series
$\wt\DY_{-N}(\gl_N)$ at the critical level $c=-N$ by
\begin{multline}\label{smuuyadou}
\wt\TT_{\mu}(u)=\tr^{}_{1,\dots,m}\ts \Ec^{}_{\Uc}\ts T^+_1(u+h\tss c_1)\dots T^+_m(u+h\tss c_m)\\[0.5em]
{}\times T_m\big(u+h\tss c_m-h\tss N/2\big)^{-1}\dots T_1\big(u+h\tss c_1-h\tss N/2\big)^{-1},
\end{multline}
where $c_a=c_a(\Uc)$ is the content
of the box occupied by $a\in\{1,\dots,m\}$ in the standard tableau $\Uc$.
By the argument of \cite[Sec.~3.4]{o:qi}, the series $\wt\TT_{\mu}(u)$
does not depend on the standard tableau $\Uc$ of shape $\mu$.

\bth\label{thm:yachardou}
All coefficients of $\wt\TT_{\mu}(u)$ belong to the center
of the $h$-adically completed algebra $\wt\DY_{-N}(\gl_N)$.
\eth

\bpf
We need to show that
\beql{comttp}
T_0(z)\ts\wt\TT_{\mu}(u)=\wt\TT_{\mu}(u)\ts T_0(z)\Fand
T^+_0(z)\ts\wt\TT_{\mu}(u)=\wt\TT_{\mu}(u)\ts T^+_0(z).
\eeq
Repeat the corresponding part of the
proof of Theorem~\ref{thm:tmuinva} and use the relations
\begin{multline}
T_0(z)\ts\overline R_{0\tss a}(z-u-h\tss c_a+h\tss N/2)
\ts T_a\big(u+h\tss c_a-h\tss N/2\big)^{-1}\\[0.4em]
{}=T_a\big(u+h\tss c_a-h\tss N/2\big)^{-1}\tss \overline R_{0\tss a}(z-u-h\tss c_a+h\tss N/2)
\ts T_0(z)
\non
\end{multline}
implied by \eqref{RTT} to get
\ben
T_0(z)\ts\wt\TT_{\mu}(u)=\wt\TT'_{\mu}(u)\ts T_0(z),
\een
where we set
\begin{multline}
\wt\TT'_{\mu}(u)=\tr^{}_{1,\dots,m}\ts \Ec^{}_{\Uc}\ts
\overline R_{0\tss 1}^{\ts-1}\dots\overline R_{0\tss m}^{\ts-1}
\ts T^+_1(u+h\tss c_1)\dots T^+_m(u+h\tss c_m)\\[0.5em]
{}\times T_m\big(u+h\tss c_m-h\tss N/2\big)^{-1}\dots T_1\big(u+h\tss c_1-h\tss N/2\big)^{-1}
\overline R_{0\tss m}\dots \overline R_{0\tss 1}.
\non
\end{multline}
The same argument as in the proof of Theorem~\ref{thm:tmuinva} shows that
$\wt\TT'_{\mu}(u)=\wt\TT_{\mu}(u)$ thus verifying the first relation in \eqref{comttp}.
A similar calculation verifies the second relation. It relies on the identity
\ben
T^+_0(z)\ts T^+_a(u+h\tss c_a)=\overline R_{0\tss a}
(z-u-h\tss c_a)^{-1}\ts T^+_a(u+h\tss c_a)\ts T^+_0(z)
\ts \overline R_{0\tss a}(z-u-h\tss c_a)
\een
implied by \eqref{RTT}, and
\begin{multline}
T^+_0(z)\ts\overline R_{0\tss a}(z-u-h\tss c_a)\ts T_a\big(u+h\tss c_a-h\tss N/2\big)^{-1}\\[0.5em]
{}=T_a\big(u+h\tss c_a-h\tss N/2\big)^{-1}\ts \overline R_{0\tss a}(z-u-h\tss c_a+h\tss N)
\ts T^+_0(z)
\non
\end{multline}
which follows from \eqref{RTTdou} with the use of
\eqref{up}.
\epf

The following formula for $\wt\TT_{\mu}(u)$ in the case where $\mu=(1^N)$
is a column diagram is a consequence of \eqref{rmatdet} and its
counterpart for the matrix $T(u)$.

\bpr\label{prop:qdethh}
At the critical level $c=-N$ we have
\ben
\wt\TT_{(1^N)}(u)=\qdet T^+(u)\big(\qdet T(u-h\tss N/2)\big)^{-1}.
\vspace{-0.6cm}
\een
\qed
\epr

By applying $\wt\TT_{\mu}(u)$ to the vacuum vector of the module
$\Vc_{-N}(\gl_N)$ we get
\beql{tvac}
\wt\TT_{\mu}(u)\tss\vac=\TT^+_{\mu}(u)\tss\vac,
\eeq
where
\beql{smuuyadualh}
\TT^+_{\mu}(u)=\tr^{}_{1,\dots,m}\ts \Ec^{}_{\Uc}\ts T^+_1(u+h\tss c_1)\dots T^+_m(u+h\tss c_m),
\eeq
in accordance with Sec.~\ref{subsec:ivmdy}. In particular, Theorem~\ref{thm:tmuinva}
follows from Theorem~\ref{thm:yachardou}; cf. \cite{fjmr:hs}.
As another application of \eqref{tvac}, we get the following.

\bco\label{cor:commut}
The coefficients of all series $\TT^+_{\mu}(u)$ generate a commutative
subalgebra of the $h$-adically completed dual Yangian $\Y^+(\gl_N)$.
\eco

\bpf
Let $\mu$ and $\nu$ be partitions having at most $N$ parts.
By Theorem~\ref{thm:yachardou} we have
\ben
\wt\TT_{\mu}(u)\wt\TT_{\nu}(u)\tss\vac
=\wt\TT_{\mu}(u)\TT^+_{\nu}(u)\tss\vac
=\TT^+_{\nu}(u)\tss \wt\TT_{\mu}(u)
\tss\vac=\TT^+_{\nu}(u)\tss \TT^+_{\mu}(u)\tss\vac.
\een
Swapping the operators, we conclude
that the coefficients of the series $\TT^+_{\mu}(u)$ and $\TT^+_{\nu}(u)$
pairwise commute in the dual Yangian.
\epf

\bre\label{rem:quace}
By the definition \eqref{qva1} of the vertex operators, evaluating \eqref{smuuyadualh}
at $u=0$ we get
$Y\big(\mathbb{T}^{+}_{\mu} (0)\vac,z\big)=\wt{\mathbb{T}}_{\mu} (z)$,
where $\wt{\mathbb{T}}_{\mu} (z)$ is given by \eqref{smuuyadou},
but the coefficients of this series are now understood as operators
on the vacuum module; cf. \cite[Sec.~3.2.2]{f:lc}.
\qed
\ere

\subsection{Center of the quantum affine vertex algebra}
\label{subsec:cava}

By Proposition~\ref{prop:associativealgebra}, the center of a quantum
vertex algebra is an associative algebra with respect to the product defined
in \eqref{as:prod}.
Moreover,
due to Proposition~\ref{prop:scommut},
this algebra is
$\mathcal{S}$-commutative, i.e. its elements satisfy  \eqref{sloc}.
The results of this section will imply that the center of the
quantum affine vertex algebra $\Vc_c(\gl_N)$
associated with $\gl_N$ is {\em commutative}, so it
shares the commutativity property of the center of
a vertex algebra; cf. \cite[Lemma 3.3.2]{f:lc}.
It follows from the definition \eqref{center} that the center
coincides with the subspace of invariants
\beql{ffvacyh}
\z\big(\Vc_c(\gl_N)\big)=
\{v\in \Vc_c(\gl_N)\ |\ t_{ij}^{\ts(r)}\tss v=0\qquad
\text{for $r\geqslant 1$\ \  and all\ \  $i,j$}\}
\eeq
of the vacuum module $\Vc_c(\gl_N)$;
cf. \eqref{ffvacy}. Hence, $\z\big(\Vc_c(\gl_N)\big)$ can be identified
with a subspace of the $h$-adically completed dual Yangian $\Y^+(\gl_N)$. Moreover, it follows from
\eqref{qva1} that
the product \eqref{as:prod} on the center coincides with the product
in the algebra $\Y^+(\gl_N)$. Therefore, by Proposition~\ref{prop:associativealgebra}
the center can be regarded as a $D$-invariant
associative subalgebra of the dual Yangian.

Now assume that the level is critical,
$c=-N$, and set $\Vc_{\text{cri}}=\Vc_{-N}(\gl_N)$.
In Corollaries~\ref{cor:symasym} and \ref{cor:newton} we constructed
three families of invariants of the extended vacuum module at the critical level.
In accordance with the definition \eqref{center}, we can reformulate
these results for the current setting by stating that
all coefficients of the series
\ben
\tr^{}_{1,\dots,m}\ts H^{(m)}\ts T^+_1(u-h\tss m+h)\dots T^+_m(u),\qquad
\tr^{}_{1,\dots,m}\ts A^{(m)}\ts T^+_1(u)\dots T^+_m(u-h\tss m+h)
\een
and
\ben
\tr\ts T^+(u)\dots T^+(u-h\tss m+h),
\een
belong to the center $\z(\Vc_{\text{cri}})$ of the quantum affine
vertex algebra $\Vc_{\text{cri}}$ (we have used the shift $u\mapsto u-h\tss m+h$
for the first series).
We will use these families to produce
generators of $\z(\Vc_{\text{cri}})$. Extend $\Vc_{\text{cri}}$
to a module over the field $\CC((h))$ and introduce its
elements as coefficients of the series
\ben
\bal
\Phi_m(u)&=h^{-m}\ts \sum_{k=0}^m(-1)^k\ts \binom{N-k}{m-k}\ts
\tr^{}_{1,\dots,k}\ts A^{(k)}\ts T^+_1(u)\dots T^+_k(u-h\tss k+h),\\
\Psi_m(u)&=h^{-m}\ts \sum_{k=0}^m(-1)^k\ts \binom{N+m-1}{m-k}\ts
\tr^{}_{1,\dots,k}\ts H^{(k)}\ts T^+_1(u-h\tss k+h)\dots T^+_k(u),
\eal
\een
and
\ben
\Theta_m(u)=h^{-m}\ts \sum_{k=0}^m(-1)^k\ts\binom{m}{k}\ts
\tr\ts T^+(u)\dots T^+(u-h\tss k+h).
\een
Define the coefficients by
\ben
\Phi_m(u)=\sum_{r=0}^{\infty}\ts \Phi^{(r)}_{m}\tss u^r,\qquad
\Psi^{}_{m}(u)=\sum_{r=0}^{\infty}\ts \Psi^{(r)}_{m}\tss u^r
\Fand
\Theta^{}_{m}(u)=\sum_{r=0}^{\infty}\ts \Theta^{(r)}_{m}\tss u^r.
\een

\bpr\label{prop:phipsithe}
All coefficients of the series $\Phi_m(u)$, $\Psi^{}_{m}(u)$ and $\Theta^{}_{m}(u)$
belong to the $\CC[[h]]$-module $\z(\Vc_{\text{\rm cri}})$.
Moreover, each family $\Phi^{(r)}_{m}$, $\Psi^{(r)}_{m}$
and $\Theta^{(r)}_{m}$
with $m=1,\dots,N$ and $r=0,1,\dots$ is algebraically independent.
\epr

\bpf
First consider the series $\Phi_m(u)$. As in Sec.~\ref{subsec:ivmdy} embed
the dual Yangian into the algebra of formal series $\Y^+(\gl_N)[[u,\di_u]]$
and introduce the element
\beql{onemtdh}
\tr^{}_{1,\dots,m}\tss A^{(m)} \big(1-T^+_1(u)\tss e^{-h\tss \di_u}\big)\dots
\big(1-T^+_m(u)\tss e^{-h\tss \di_u}\big)
\eeq
as in \eqref{onemtd}. By repeating the corresponding argument
in Sec.~\ref{subsec:ivmdy} we find
that the element \eqref{onemtdh} coincides with
\beql{parteqh}
\sum_{k=0}^m(-1)^k \binom{N-k}{m-k}\ts
\tr^{}_{1,\dots,k}\ts A^{(k)}\ts
T^+_{1}(u)\dots T^+_{k}(u-h\tss k+h)\tss e^{-k\tss h\tss\di_u}.
\eeq
Observe that the constant term of \eqref{parteqh}, as a formal power series in $\di_u$,
coincides with $h^m\tss \Phi_m(u)$.
On the other hand, each factor in
\eqref{onemtdh} takes the form
\ben
1-T^+_i(u)\tss e^{-h\tss \di_u}\equiv h\ts (\di_u+T^{(-1)}+T^{(-2)}\tss u+\dots)
\mod h^2\ts \Vc_{\text{cri}},
\een
where $T^{(-r)}=[\tss t_{ij}^{(-r)}]$ is the matrix of generators.
This shows that the series $h^m\tss \Phi_m(u)$ belongs to $h^m\tss \Vc_{\text{\rm cri}}[[u]]$
and so all coefficients of $\Phi_m(u)$ belong to the
$\CC[[h]]$-module $\z(\Vc_{\text{cri}})$. Furthermore, taking the classical limit
$h=0$ we find that the image of the series $T^{(-1)}+T^{(-2)}\tss u+\dots$
in the algebra $\U\big(t^{-1}\gl_N[t^{-1}]\big)[[u]]$ coincides
with $E_+(u)$ as defined in \eqref{epu}. By Corollary~\ref{cor:ffah},
the family of elements $\phi^{(r)}_{m\tss m}$, found as constant terms
of the polynomials \eqref{parteq} in $\di_u$, is algebraically independent.
Hence so is the family of the coefficients $\Phi^{(r)}_{m}$.
Indeed, if there is a polynomial with coefficients in $\CC[[h]]$ providing
an algebraic dependence of the $\Phi^{(r)}_{m}$, we may assume that at least one
of its coefficients is not zero modulo $h$. Then the evaluation $h=0$
makes a contradiction.

The arguments for the families $\Psi^{}_{m}(u)$ and $\Theta^{}_{m}(u)$
are quite similar. One additional observation for the family $\Psi^{}_{m}(u)$
is the identity
\ben
\tr^{}_{1,\dots,m}\ts H^{(m)}\ts T^+_1(u-h\tss m+h)\dots T^+_m(u)
=\tr^{}_{1,\dots,m}\ts\ts T^+_1(u)\dots T^+_m(u-h\tss m+h)\tss H^{(m)}.
\een
It follows by applying the fusion formula \eqref{fusion}
for $H^{(m)}$, then the defining relations \eqref{RTTdual}
and the conjugation by the longest permutation of $\Sym_m$.
\epf

We can now prove a quantum analogue of
the Feigin--Frenkel theorem~\cite{ff:ak}; see Sec.~\ref{subsec:ivmdy}.

\bth\label{thm:qff}
The center at the critical level $\z(\Vc_{\text{\rm cri}})$
is a commutative algebra. It is topologically generated by each
of the families $\Phi^{(r)}_{m}$, $\Psi^{(r)}_{m}$
and $\Theta^{(r)}_{m}$
with $m=1,\dots,N$ and $r=0,1,\dots$
\eth

\bpf
First we point out that the coefficients of all series
$\Phi_m(u)$, $\Psi^{}_{m}(u)$ and $\Theta^{}_{m}(u)$ pairwise commute.
This is well-known for the Yangian counterparts of the series introduced
in Corollaries~\ref{cor:symasym} and \ref{cor:newton}
(with the matrix $T^+(u)$ replaced with $T(u)$)
in relation with Bethe subalgebras \cite{ks:qs}; see also \cite[Ch.~1]{m:yc}.
The same proof applies for the dual Yangian.
Alternatively, this fact is obtained as a consequence of Corollary~\ref{cor:commut}.

Now suppose that $w\in \z(\Vc_{\text{\rm cri}})$.
We will prove by induction that for all $n\geqslant 0$ there exists a polynomial $Q$ in the variables
$\Phi^{(r)}_{m}$ such that
$
w-Q\in h^n \Vc_{\text{\rm cri}}.
$
Assuming that this holds for some $n\geqslant 0$, write
\ben
w-Q=h^n\tss w_n+h^{n+1}\tss w_{n+1}+\dots \qquad  \text{with }w_{k} \in V_0,
\een
where we assume that $\Vc_{\text{\rm cri}}=V_0[[h]]$.
Since $w-Q$ belongs to the center of the vacuum module, we can conclude that
$w_n\in \z(\Vc_{\text{\rm cri}})\mod h$. Taking the classical limit $h=0$ we find
that the image $\overline w_n$ of $w_n$ in $\U\big(t^{-1}\gl_N[t^{-1}]\big)$
belongs to the Feigin--Frenkel center $\z(\wh\gl_N)$. Therefore,
$\overline w_n$ is a polynomial $S$ in the variables $\phi^{(r)}_{m\tss m}$;
see Corollary~\ref{cor:ffah}. Replace these variables with the respective
elements $\Phi^{(r)}_{m}$ to get a polynomial $S'\in \z(\Vc_{\text{\rm cri}})$.
The difference $w_n-S'$ belongs to $h\tss \Vc_{\text{\rm cri}}$.
Therefore,
\ben
w-Q-h^n\tss S'\in h^{n+1} \Vc_{\text{\rm cri}},
\een
which completes the induction argument.

Thus, any element $w\in \z(\Vc_{\text{\rm cri}})$ can be approximated
by polynomials in the variables $\Phi^{(r)}_{m}$ and so they are topological
generators of the center.
The same argument works for the other two
families. In particular, this implies that the algebra
$\z(\Vc_{\text{\rm cri}})$ is commutative.
\epf

Finally, consider the quantum affine vertex algebra $\Vc_c(\gl_N)$ with
$c\ne -N$.
The center of the affine
vertex algebra $V_{\ka}(\gl_N)$ with $\ka\ne -N$
is known to be generated by the elements
\beql{sdegone}
E_{11}[-r-1]+\dots+E_{NN}[-r-1],\qquad r=0,1,\dots.
\eeq
By Proposition~\ref{prop:qdet},
the coefficients of the quantum determinant
\beql{qdeth}
\qdet T^+(u)=\sum_{\si\in\Sym_N}\sgn \si\cdot
t^+_{ \si(1)\tss 1}(u)\dots t^+_{ \si(N)\tss N}(u-h\tss N+h),
\eeq
as defined in \eqref{qdetplus},
belong to the center $\z\big(\Vc_c(\gl_N)\big)$. Write
\ben
\qdet T^+(u)=1-h\tss\big(d_0+d_1\tss u+d_2\tss u^2+\dots\big).
\een
Under the
classical limit $h\to 0$, the image of $d_r$ in $\U\big(t^{-1}\gl_N[t^{-1}]\big)$
coincides with the element \eqref{sdegone}. The same argument
as in the proof of Theorem~\ref{thm:qff} yields the following.

\bpr\label{prop:noncr}
The center $\z\big(\Vc_c(\gl_N)\big)$ with $c\ne -N$
is a commutative algebra. It is topologically generated by
the family $d_0,d_1,\dots$ of algebraically independent elements.
\qed
\epr

In particular, $\z\big(\Vc_c(\gl_N)\big)$ is a commutative algebra
for all values of $c$. It seems to be plausible that this commutativity property of the center extends
to all quantum VOAs; cf. Remark \ref{remark-new}.


\newpage

\small

\noindent
N.J.:\qquad\qquad\qquad\qquad\\
Department of Mathematics\\
North Carolina State University, Raleigh, NC 27695, USA\\
jing@math.ncsu.edu\\

\noindent
N.J. \& F.Y.:\newline
School of Mathematical Sciences\\
South China University of Technology\\
Guangzhou, Guangdong 510640, China\\


\noindent
S.K. \& A.M.:\newline
School of Mathematics and Statistics\newline
University of Sydney,
NSW 2006, Australia\newline
kslaven@maths.usyd.edu.au\newline
alexander.molev@sydney.edu.au\\

\noindent
S.K.:\\
Department of Mathematics\\
University of Zagreb, 10000 Zagreb, Croatia

\end{document}